\documentclass[reqno]{amsart}
\usepackage{amsaddr}
\usepackage[pagewise]{lineno}
\usepackage{url}
\usepackage{xcolor}
\usepackage[font=scriptsize,labelfont=bf]{caption}
\usepackage[symbol]{footmisc}
\usepackage[latin1]{inputenc}
\usepackage{amsmath}
\usepackage{xcolor}
\usepackage{amssymb}
\usepackage{enumerate}
\usepackage{graphicx}
\usepackage{amsthm}
\usepackage{amsopn}
\usepackage{amsfonts}
\usepackage{upgreek}
\usepackage{cite}
\usepackage{amsfonts}
\usepackage{amssymb}
\usepackage{caption}
\usepackage{amsthm}
\usepackage{mathtools, amssymb}
\usepackage[all]{xy}
\numberwithin{equation}{section}

\newcommand{\U}{\mathbb U}

\providecommand{\U}[1]{\protect\rule{.1in}{.1in}}

\newtheorem{lemma}{Lemma}
\newtheorem{theorem}{Theorem}
\newtheorem{definition}{Definition}
\newtheorem{remark}{Remark}

\begin{document}
\title[A variation of parameters formula for  nonautonomous IDEPCAG]{A variation of parameters formula for  nonautonomous linear impulsive differential equations with piecewise constant arguments of generalized type
}
\author{Ricardo Torres*}
\address{
 Instituto de Cs. F\'isicas y Matem\'aticas, Facultad de Ciencias \\
Universidad Austral de Chile\\
Campus Isla Teja, Valdivia, Chile \vspace{0.3cm}\\
 Instituto de Ciencias, \\
Universidad Nacional de Gral. Sarmiento \\
Los Polvorines, Bs. Aires, Argentina}
\email{ricardo.torres@uach.cl (Corresponding author)}

\author{Manuel Pinto}
\address{
 Departamento de Matem\'aticas, Facultad de Ciencias \\
Universidad de Chile\\
Campus Las Palmeras, Santiago, Chile}
\email{pintoj.uchile@gmail.com}

\begin{abstract}
In this work, we give a variation of parameters formula for nonautonomous linear impulsive differential equations with piecewise constant arguments of generalized type. We cover several cases of differential equations with deviated arguments investigated before as particular cases. We also give some examples showing the applicability of our results. 
\end{abstract}
\maketitle

\section{Introduction}
\footnote[2]{This manuscript is dedicated to the memory of Prof. Nicol\'as Yus Su\'arez.}Occasionally, natural phenomena must be modeled using differential equations that may have discontinuous solutions, such as a piecewise constant, or the impulsive effect must be present. Some examples of such modeling can be found in the works of S. Busenberg and K. Cooke \cite{COOKE1} (where the authors modeled vertical transmission diseases) and L. Dai and M.C. Singh \cite{DAI_SINGH} (oscillatory motion of spring-mass systems subject to piecewise constant forces such $Ax([t])$ or $A\cos([t])).$ The last work studied the motion of mechanisms modeled by $$mx''(t)+kx_1=A sin\left(\omega \left[\dfrac{t}{T}\right]\right),$$
where $[\cdot]$ is the greatest integer function. 
 (See \cite{DAI}).\\
 
In the 70's, A. Myshkis \cite{203} studied differential equations with deviating arguments ($h(t)\leq t$, such as $h(t)=[t]$ or $h(t)=[t-1])$.  The Ukrainian mathematician M. Akhmet generalized those systems, introducing  differential equations of the form
\begin{equation}
y^{\prime}(t)=f(t,y(t),y(\gamma(t))),\label{depcag_eq}
\end{equation}
where $\gamma(t)$ is a \emph{piecewise constant argument of generalized type}.
In order to define such $\gamma$, let $\left(t_{n}\right)_{n\in\mathbb{Z}}$ and $\left(\zeta_{n}\right)_{n\in\mathbb{Z}}$ such that $t_{n}<t_{n+1}\, , \forall n\in\mathbb{Z}$ with $\displaystyle{\lim_{n\rightarrow\infty}t_{n}=\infty}$, $\displaystyle{\lim_{n\rightarrow -\infty}t_{n}=-\infty}$ and $\zeta_{n}\in[t_{n},t_{n+1}].$
Then,  $\gamma(t)=\zeta_{n},$ if $t\in I_{n}=\left[t_{n},t_{n+1}\right).$
I.e., $\gamma(t)$ is a step function. An elementary example of such functions is $\gamma(t) =[t]$  which is constant in every interval $[n,n+1[$ with  $n\in \mathbb{Z}$ (see \eqref{idepca_parte_entera}).\\
If a piecewise constant argument is used, the interval $I_n$ is decomposed into an advanced and delayed subintervals $I_{n}=I_{n}^{+}\bigcup I_{n}^{-}$, where
$I_{n}^{+}=[t_{n},\zeta_{n}]$ and $I_{n}^{-}=[\zeta_{n},t_{n+1}].$
This class of differential equations is known as \emph{Differential Equations with Piecewise Constant Argument of Generalized Type} (\emph{DEPCAG}). They have continuous solutions, even though $\gamma$ is discontinuous. If we assume continuity of the solutions of \eqref{depcag_eq}, integrating from $t_n$ to $t_{n+1}$, we define a finite-difference equation, so we are in the presence of a hybrid dynamic (see \cite{AK2, P2011}).\\
For example, taking $\gamma(t) =\left[ \frac{t+l}{h}\right]h $ with $0\leq l<h$, we have
\begin{eqnarray*}
\left[ \frac{t+l}{h}\right]h=nh, \text{ when } t\in I_n=[nh-l,\left(n+1\right)h-l).
\end{eqnarray*}
Then, we see that  $\gamma (t) -t \geq 0\text{ }\Leftrightarrow \text{ }t\leq nh$ and
$\gamma (t) -t \leq 0\text{ }\Leftrightarrow \text{ }t\geq nh$. Hence, we have
\begin{equation*}
I_{n}^{+}=[nh-l,nh ],\quad I_{n}^{-}=[nh ,\left(n+1\right)h-l].
\end{equation*}
 Now, if an impulsive condition is defined at $\{t_n\}_{n\in\mathbb{Z}}$, we are in the presence of the 
\emph{Impulsive differential equations with piecewise constant argument of generalized type} (\emph{IDEPCAG}) (see \cite{AK3}),
\begin{align}
x^{\prime}(t)&=f(t,x(t),x(\gamma(t))),\qquad t\neq t_{n} \nonumber \\
\Delta x(t_{n})&:=x(t_{n})-x(t_{n}^{-})=J_{n}(x(t_{n}^{-})),\qquad t=t_{n},\quad n\in\mathbb{N} \label{idepcag_gral},
\end{align}
where $x(t_n^-)=\displaystyle{\lim_{t\to t_n^-}x(t),}$ and $J_n$ is the impulsive operator (see \cite{Samoilenko}).\\

When the piecewise constant argument used in a differential equation is explicit, it will be called DEPCA (IDEPCA if it has impulses). 

\bigskip
\subsection*{An elementary and illustrative example of IDEPCA}~\par
\vskip1mm
Consider the scalar IDEPCA
\begin{align}
x'(t)&=(\alpha-1)x([t]),\qquad t\neq n \nonumber \\
x(n)&=\beta x(n^{-}),\qquad t=n,\quad n\in\mathbb{N} \label{idepca_parte_entera}.
\end{align}
where $\alpha,\beta \in\mathbb{R}, \, \beta\neq 1.$\\ 
If $t\in[n,n+1)$ for some $n\in\mathbb{Z}$, equation \eqref{idepca_parte_entera} can be written as
\begin{equation}
    x'(t)=(\alpha-1)x(n). \label{idepca_parte_entera_1}
\end{equation}
In the following, we will assume $t_0=0$. Now, integrating on $[n,n+1)$ from $n$ to $t$ we see that
\begin{equation}
    x(t)=x(n)(1+(\alpha-1)(t-n)). \label{idepca_parte_entera_2}
\end{equation}
Next, assuming continuity at $t=n+1$, we have
\begin{equation*}
    x((n+1)^{-})=\alpha x(n). \label{idepca_parte_entera_3}
\end{equation*}
Applying the impulsive condition to the last expression, we get the following \emph{finite-difference equation}
\begin{equation*}
    x((n+1))=(\alpha \beta)x(n). \label{idepca_parte_entera_3_2}
\end{equation*}
Its solution is
\begin{equation}
    x(n)=(\alpha \beta)^{n}x(0). \label{idepca_parte_entera_4}
\end{equation}
Finally, applying \eqref{idepca_parte_entera_4} in \eqref{idepca_parte_entera_2} we have
\begin{equation}
    x(t)=\left(\alpha \beta\right)^{[t]}(1+(\alpha-1)(t-[t]))x(0). \label{idepca_parte_entera_5}
\end{equation}
\begin{remark}
\begin{enumerate}
    \item[] 
    \item From \eqref{idepca_parte_entera_5}, we can conclude that the underlying dynamic is of mixed type. The discrete and the continuous parts of the system are dependent. For example, A stable continuous part (associated with the coefficient $\alpha$) can be unstabilized by the discrete part (associated with the parameter $\beta$). See \cite{Samoilenko}. 
\end{enumerate}

\end{remark}
In the next table, we describe some of the behavior of the solutions of \eqref{idepca_parte_entera_5}:

\begin{table}[h]
\begin{tabular}{|ll|}
\hline
\textbf{Behavior of solutions}  & \textbf{Condition}\\\hline
$|x(t)| \xrightarrow{t\to\infty} 0$ exponentially.  &  $|\alpha \beta|<1$ and $\alpha \beta\neq 0$.  \\
$x(t)$ is constant. & $\alpha \beta=0$  or $\alpha=\beta=1$  \\
$x(t)$ is oscillatory. & $\alpha \beta<0$   \\
$x(t)$ is piecewise constant. & $\alpha =1$\\
$|x(t)|$ is piecewise constant and  $x(t)\xrightarrow{t\to\infty} +\infty.$ & $\alpha=1$ and $|\beta|>1$\\
$x(t)$ is piecewise constant and  $x(t) \xrightarrow{t\to\infty} 0$.  & $\alpha=1$ and $0<\beta<1$\\
$|x(t)| \xrightarrow{t\to\infty} +\infty$ exponentially. & $|\alpha \beta|>1$.\\
\hline
\end{tabular}
\caption{Behavior of solutions of \eqref{idepca_parte_entera_5}}
\end{table}
\begin{figure}[h!]
\centering
\includegraphics[scale=0.3]{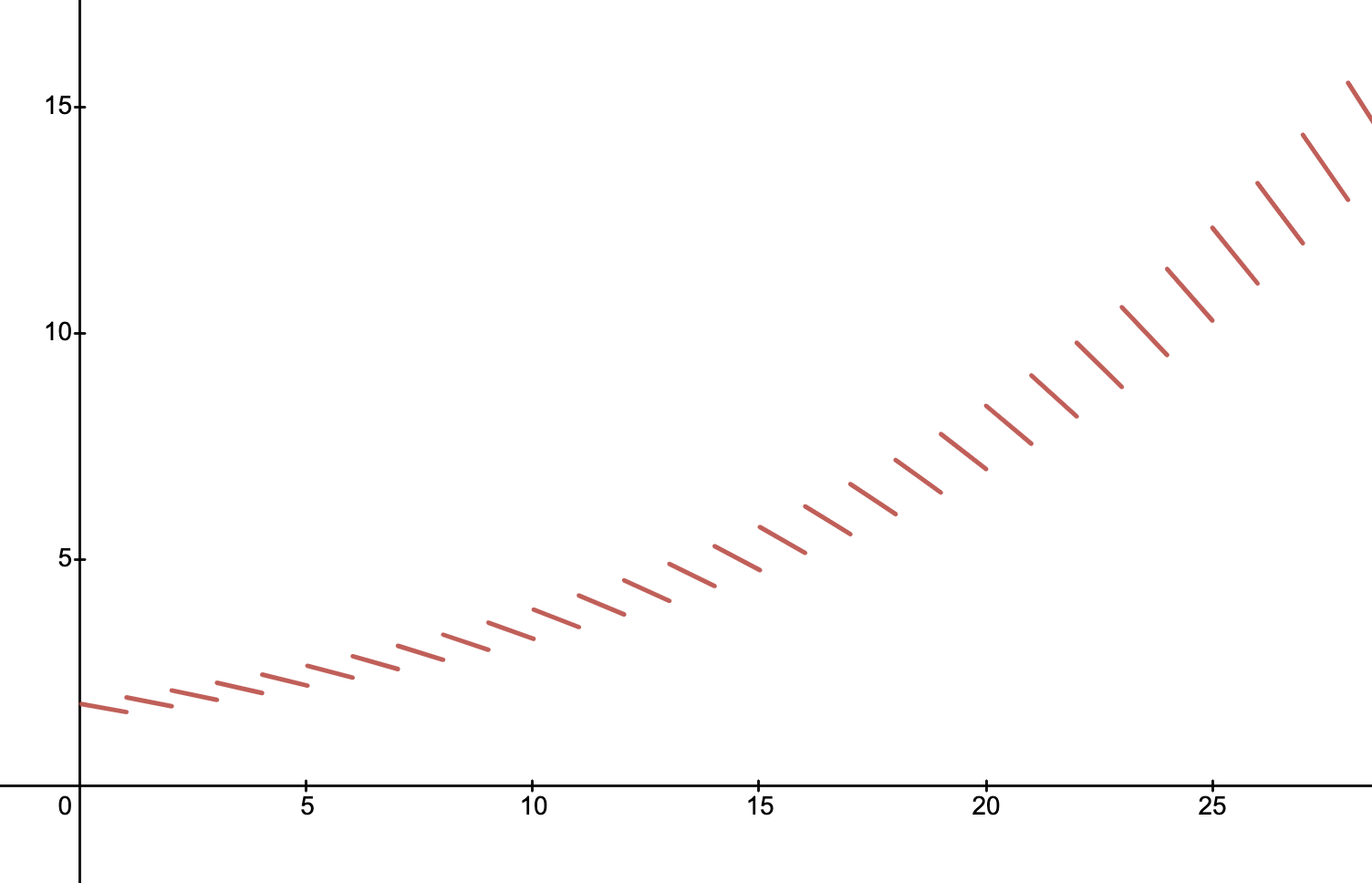}
\caption{Solution of \eqref{idepca_parte_entera} with
$\alpha=0.9 $, $\beta=1.2$, $x_0=1.8$.}
\end{figure}
\begin{figure}[h!]
\centering
\includegraphics[scale=0.3]{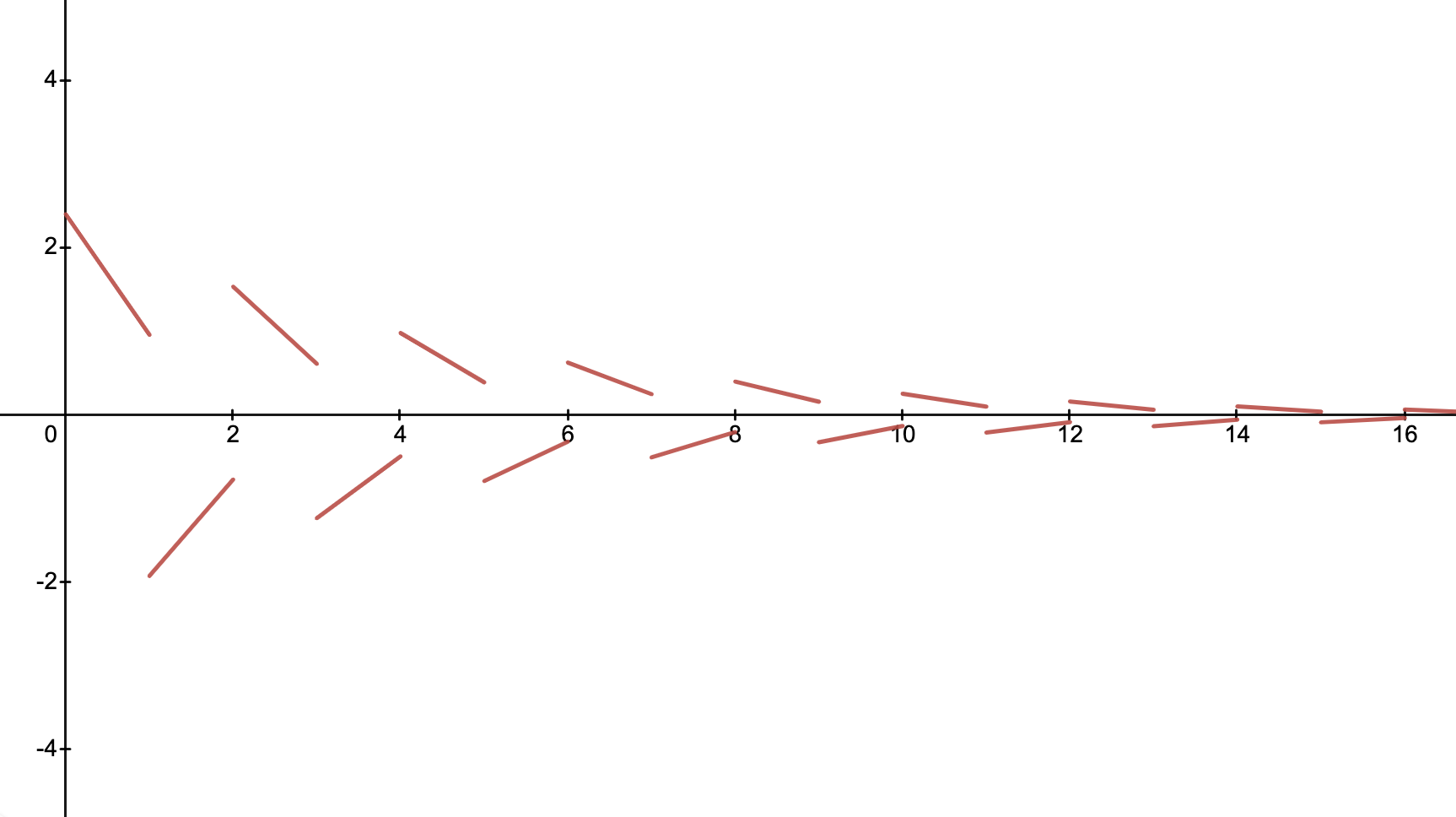}
\caption{solution of \eqref{idepca_parte_entera} with
$\alpha=0.4 $, $\beta=-2$, $x_0=2.4$.}
\end{figure}

\newpage
\subsection{Why study IDEPCAG?: impulses in action}\label{ejemplo_berek_avanzado}~\par
\vskip1mm
\textbf{Example 1.}\,\,Let the following scalar linear DEPCA
\begin{equation}
    x'(t)=a(t)(x(t)-x([t]),\quad x(\tau)=x_0,\label{berek_depca}
\end{equation}
and the scalar linear IDEPCA
\begin{equation}
    \begin{tabular}{ll}
$z'(t)=a(t)\left(z(t)-z([t])\right),$ & $t\neq k$ \\ 
$z(k)=c_k z(k^{-}),$ & $t=k,\quad k\in\mathbb{Z},$ 
\end{tabular}
\label{berek_sistema_avanzado}
\end{equation}
where $a(t)$ is a continuous locally integrable function and $(c_k)_{k\in\mathbb{N}}$ a real sequence such that $c_k\notin \{0,1\}$, for all $k\in\mathbb{N}$. As $\gamma(t)=[t],$ we have $t_{k}=k=\zeta_k=k$ if $t\in[k,k+1),\,k\in\mathbb{Z}$.\\
The solution of \eqref{berek_depca} is $x(t)=x_0,\,\, \forall t\geq \tau.$ 
I.e., all the solutions are constant (see \cite{P2011}). \\
On the other hand, as we will see, the solution of \eqref{berek_sistema_avanzado} is
$$z(t)=\left(\prod_{j=k(\tau)+1}^{k(t)}c_j\right)z(\tau),\quad t\geq \tau,$$
where $k(t)=k$ is the only integer such that $t\in[k,k+1]$.\\
Hence, all the solutions are nonconstant if $c_j\neq 1$ and $c_j\neq 0$, for all $j\geq k(\tau)$. This example shows the differences between DEPCA and IDEPCA systems. The discrete part of the system can greatly impact the whole dynamic,  determining the qualitative properties of the solutions.
\begin{figure}[h!]
\centering
\includegraphics[scale=0.3]{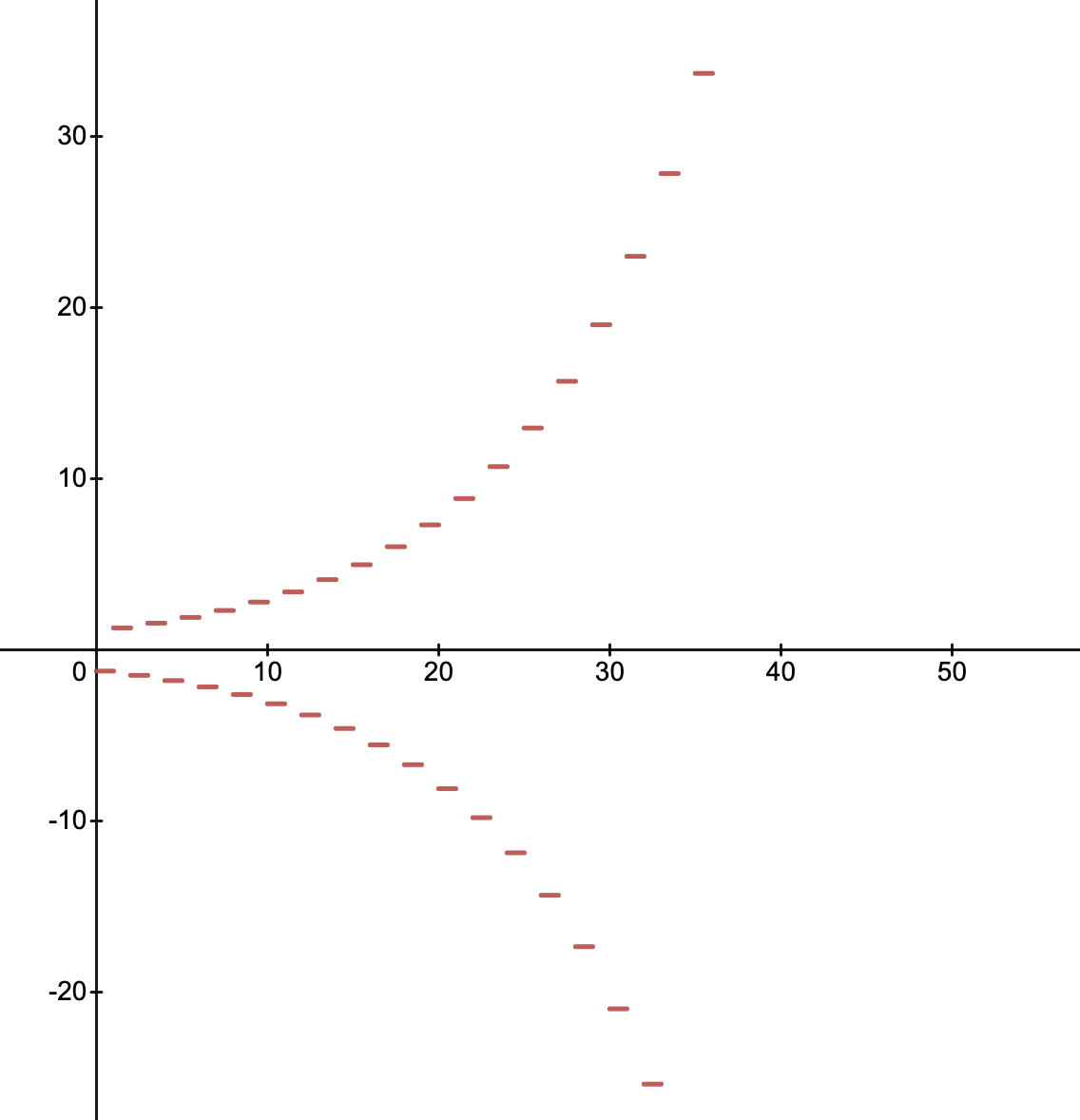}
\caption{Solution of \eqref{berek_sistema_avanzado} with $c_k=-1.1$ and $z(0)=-1.2$}
\end{figure}

\subsection{Fundamental matrices and variation of parameters formulas: an overview}\label{section_overview}~\par
\vskip1mm
\subsubsection{The fundamental matrix of a DEPCA system}~\par
\vskip1mm
In \cite{CoWi84}, K.L. Cooke and J. Wiener were the first to obtain a fundamental matrix for a scalar \emph{DEPCA}'s using the delayed piecewise constant arguments $\gamma(t) =[t]$, $\gamma(t)=[t-1]$, $\gamma (t)=[t-n]$ and $\gamma (t)=t-n[t].$ Also, they considered the very interesting scalar DEPCA
\begin{align*}
x^{\prime }(t)=a(t)x(t)+\displaystyle{\sum_{i=0}^{n}a_i(t)x([t-i])},\quad a_n\neq 0,
\end{align*}
and 
\begin{align*}
x^{\prime }(t)=ax(t)+\displaystyle{\sum_{i=1}^{n}a_ix(t-i[t])}.
\end{align*}
Also, in \cite{ShWi83}$,$ K.L. Cooke and S.M. Shah studied the DEPCA
\begin{align*}
x^{\prime }(t)=a(t)x(t)+\displaystyle{\sum_{i=0}^{n}a_k(t)x([t+i])},\quad n\leq 2.
\end{align*}
Then, in \cite{CoWi87}, K.L. Cooke and J. Wiener studied the mixed-type piecewise constant argument  
$\gamma (t) =2\left[ \frac{t+1}{2}\right] $ and considered the DEPCA
\begin{align*}
z^{\prime }(t)& =az(t)+bz(2[\left( t+1\right) /2])
\end{align*}
Additionally, in \cite{WiAf88}, K.L. Cooke and A.R. Aftabizadeh  considered the mixed-type piecewise constant argument  
$\gamma (t) =m\left[ \frac{t+k}{m}\right] $ where $0<k<m$, $k,m,n\in \mathbb{Z}^{+}$, and they studied the DEPCA 
\begin{align*}
w^{\prime }(t)& =aw(t)+bw(m[\left(t+k\right) /m]).
\end{align*}~\par
\vskip1mm
\subsubsection{Variation of parameters formula for a DEPCA}~\par
\vskip1mm
In \cite{JaDeo} (1991), N. Jayasree and S.G. Deo were the first to consider the advanced and delayed parts of the solutions studying the equation
\begin{align*}
z^{\prime }(t)& =az(t)+bz(2[\left( t+1\right) /2])+f(t),
\end{align*}
obtaining a variation of parameters formula for this DEPCA, in terms of the homogeneous linear DEPCA associated:
\begin{align*}
    z(t)=&y(t)+\sum_{j=0}^{[(t+1)/2]-1}\lambda^{-1}(1)\int_{2j}^{2j+1}\Psi(t,2j)\phi(2j+1,s)f(s)ds\\
    &-\sum_{j=1}^{[(t+1)/2]}\lambda^{-1}(1)\int_{2j}^{2j-1}\Psi(t,2j)\phi(2j-1,s)f(s)ds\\
    &+\int_{2[(t+1)/2]}^t\phi(t,s)f(s)ds,
\end{align*}
where 
\begin{eqnarray*}
    \lambda(t)=exp(at)\left(1+a^{-1}b\right)-a^{-1}b,
\end{eqnarray*}
$\phi$ and $\Psi$ are the fundamental solutions of
$x'(t)=ax(t)$ and $y^{\prime }(t) =ay(t)+by(2[\left( t+1\right) /2])$
respectively.\\

In \cite{MenYan} (2001), Q. Meng and J. Yan obtained a variation of parameters formula for the differential equation
\begin{equation*}
x'(t)+a(t)x(t)+b(t)x(g(t))=f (t), \text{ for } t>0
\end{equation*}
where $a(t),b(t)$ and $f(t)$ are locally integrable functions on $[0,\infty),$ $g(t)$ is a piecewise constant function defined by
$g(t)=np$ for $t\in[np-l,(n+1)p-l)$ with $n\in\mathbb{N}$ and $p,l$ positive constants such that $p>l$. The authors studied the oscillation and asymptotic stability properties of the solutions.\\

In \cite{A04} (2008), M. Akhmet considered the DEPCAG for systems
\begin{align}
z^{\prime }(t)& =A(t)z(t)+B(t)z(\gamma(t))+F(t), \label{DEPCAG_no_homogeneo} \\
w^{\prime }(t)& =A(t)w(t)+B(t)w(\gamma(t))+g\left(t,w(t),w(\gamma(t)) \right) .
\end{align}
where $A(t),B(t)\in C(\mathbb{R})$ are $n\times n$ real valued  uniformly bounded on $\mathbb{R}$ matrices, $g(t,x,y)\in C(\mathbb{R}\times \mathbb{R}^n\times \mathbb{R}^n)$ is an $n\times 1$ Lipschitz real valued function with $g(t,0,0)=0$, $\gamma(t)$ is a piecewise constant argument of generalized type. The author found the following  variation of parameters formula
\begin{eqnarray*}
w(t)=&&\displaystyle{W(t,t_0)w_0+W(t,t_0)\int_{t_0}^{\zeta_i}X(t_0,s)g(s,w(s),w(\gamma(s)))ds}\\
&&+\displaystyle{\sum_{k=i}^{j-1}W(t,t_{k+1})\int_{\zeta_k}^{\zeta_{k+1}}X(t_{k+1},s)g(s,w(s),w(\gamma(s)))ds}\\
&&+\displaystyle{\int_{\zeta_j}^{t}X(t,s)g(s,w(s),w(\gamma(s)))ds},
\end{eqnarray*}
where $j=j(t)$ is the only $j\in\mathbb{Z}$ such that $t_{j(t)}\leq t \leq t_{j(t)+1},$ $t_k\leq \zeta_k \leq t_{k+1}$, $t_i\leq t_0\leq t_{i+1}$,  $X$ is the fundamental matrix of 
\begin{equation*}
x^{\prime }(t) =A(t)x(t),
\end{equation*}
and $W$ is the fundamental matrix of the homogeneous linear DEPCAG
\begin{equation*}
y^{\prime }(t)=A(t)y(t)+B(t)y(\gamma(t)).    
\end{equation*}
Later, in \cite{P2011} (2011), M. Pinto gave a new DEPCAG variation of parameters formula. This time, the author considered the delayed and advanced intervals defined by the general piecewise constant argument
\begin{eqnarray*}
z(t)=&&\displaystyle{W(t,t_0)z_0+\underbrace{W(t,t_0)\int_{t_0}^{\zeta_i}X(t_0,s)g(s,z(s),z(\gamma(s)))ds}_{I_k^{+}}}\\
&&+\displaystyle{\sum_{k=i+1}^{j}\underbrace{W(t,t_{k})\int_{t_k}^{\zeta_{k}}X(t_{k},s)g(s,z(s),z(\gamma(s)))ds}_{I_k^{+}}}\\
&&+\displaystyle{\sum_{k=i}^{j-1}\underbrace{W(t,t_{k+1})\int_{\zeta_k}^{t_{k+1}}X(t_{k+1},s)g(s,z(s),z(\gamma(s)))ds}_{I_k^{-}}}\\
&&+\underbrace{\displaystyle{\int_{\zeta_j}^{t}X(t,s)g(s,z(s),z(\gamma(s)))ds}}_{I_k^{-}},
\end{eqnarray*}
where $t_i \leq t_0 \leq t_{i+1}$ and $t_{j(t)}\leq t \leq t_{j(t)+1}.$ \\
In the DEPCAG theory, decomposing the interval $I_n$ into the advanced and delayed subintervals is critical. As we will see, it is necessary for the forward or backward continuation of solutions. \\

\subsubsection{Variation of parameters formula for an IDEPCA: the impulsive effect applied}~\par
\vskip1mm
For the IDEPCA case, In \cite{oztepe2012convergence} (2012), G. Oztepe and  H. Bereketoglu studied the scalar IDEPCA
\begin{align}
x^{\prime}(t)&=a(t)(x(t)-x([t+1]))+f(t),\quad x(0)=x_0,\qquad t\neq n\in\mathbb{N} \nonumber \\
\Delta x(n)&=d_n ,\qquad t=n,\quad n\in\mathbb{N} \label{oztepe_berek},
\end{align}
They proved the convergence of the solutions to a real constant when $t\to\infty$, and they showed the limit value in terms of $x_0$, using a suitable integral equation. They concluded the following expression for the solutions of \eqref{oztepe_berek}
\begin{align*}
x(t)=&exp\left(\int_{[t]}^ta(u)du\right)x([t])+\left(1-exp\left(\int_{[t]}^t a(u)du\right)\right)x([t+1])\\
&+\int_{[t]}^t exp\left(\int_s^t a(u)du\right)f(s)ds,
\end{align*}
where
\begin{align*}
x([t])=&x_0+\sum_{j=0}^{[t]-1}\left(\int_j^{j+1} exp\left(-\int_j^s a(u)du\right)f(s)ds+exp\left(-\int_j^{j+1} a(u)du\right)d_{j+1}\right).
\end{align*}
For the IDEPCA case, in \cite{Kuo2023} (2023), K-S. Chiu and I. Berna considered the following impulsive differential equation with a piecewise constant argument

\begin{align}
&y^{\prime}(t)=a(t)y(t)+b(t)y\left(p\left[\frac{t+l}{p}\right]\right),\quad y(\tau)=c_0,\qquad t\neq kp-l \nonumber \\
&\Delta y(kp-l)=d_k y({kp-l}^{-}),\qquad t=kp-l,\quad k\in\mathbb{Z} \label{kuo_homogenea},
\end{align}
and
\begin{align}
&y^{\prime}(t)=a(t)y(t)+b(t)y\left(p\left[\frac{t+l}{p}\right]\right)+f(t),\quad y(\tau)=c_0,\qquad t\neq kp-l \nonumber \\
&\Delta y(kp-l)=d_k y({kp-l}^{-}),\qquad t=kp-l,\quad k\in\mathbb{Z} \label{kuo_no_homogenea},
\end{align}
where $a(t)\neq 0$, $b(t)$ and $f(t)$ are real-valued continuous functions, $p<l$ and $d_k\in \mathbb{R}-\{1\}.$ 
The authors obtained criteria for the existence and uniqueness, a variation of parameters formula, a Gronwall-Bellman inequality, stability and oscillation criteria for solutions for \eqref{kuo_homogenea} and \eqref{kuo_no_homogenea}.\\

To our knowledge, there is no variation formula for impulsive differential equations with a generalized constant argument. 
As we have shown, some authors have studied just some particular cases before.

\section{Aim of the work}
We will get a variation of parameters formula associated with \emph{IDEPCAG} system
\begin{equation}
\begin{tabular}{ll}
$x^{\prime }(t)=A(t)x(t)+B(t)x(\gamma (t))+F(t),$ & $t\neq t_{k}$ \\ 
$\Delta x|_{t=t_{k}}=C_{k}x(t_{k}^{-})+D_{k},$ & $t=t_{k},$%
\end{tabular}
\label{sistema_idepcag_general_abstract}
\end{equation}

extending the particular case treated in \cite{Kuo2023} and the general results of the DEPCAG case studied in \cite{P2011} to the IDEPCAG context.

\section{Preliminaires}
Let $\mathcal{PC}(X, Y)$ be the set of all functions $r: X\to Y$ which are continuous for $t\neq t_k$ and continuous from the left with discontinuities of the first kind at $t=t_k$. Similarly, let $\mathcal{PC}^{1}(X,Y)$ the set of functions $s:X\to Y$ such that $s'\in \mathcal{PC}(X,Y).$
\begin{definition}[DEPCAG solution]
A continuous function $x(t)$ is a solution of \eqref{depcag_eq} if:
\begin{itemize}
\item[(i)] $x'(t)$ exists at each point $t\in \mathbb{R}$ with the possible exception at the times $t_{k}$, $k\in \mathbb{Z}$, where the one side derivative exists.
\item[(ii)] $x(t)$ satisfies \eqref{depcag_eq} on the intervals of the form $(t_{k},t_{k+1})$, and it holds for the right
derivative of $x(t)$ at $t_{k}$.
\end{itemize}
\end{definition}
\begin{definition}[IDEPCAG solution]
A piecewise continuous function $y(t)$ is a solution of \eqref{idepcag_gral} if:
\begin{itemize}
\item[(i)] $y(t)$ is continuous on $I_k=[t_k,t_{k+1})$ with first kind discontinuities at $t_k, \,k\in \mathbb{Z}$, 
where $y'(t)$ exists at each $t\in \mathbb{R}$ with the possible exception at the times $t_{k}$, where lateral derivatives exist (i.e. $y(t)\in \mathcal{PC}^{1}([t_k,t_{k+1}),R^n)$).
\item[(ii)] The ordinary differential equation
$$y^{\prime}(t)=f(t,y(t),y(\zeta_{k}))$$ holds on every interval $I_{k}$, where $\gamma(t)=\zeta_k$.
\item[(iii)]
For $t=t_{k}$, the impulsive condition 
$$\Delta y(t_{k})=y(t_{k})-y(t_{k}^{-})=J_{k}(y(t_{k}^{-}))$$
holds. I.e., $y(t_{k})=y(t_{k}^{-})+J_{k}(y(t_{k}^{-}))$, where $y(t_{k}^{-})$ denotes the left-hand limit
of the function $y$ at $t_k$.
\end{itemize}
\end{definition}
Let the IDEPCAG system:
\begin{equation}
\begin{tabular}{ll}
$x^{\prime }(t)=f(t,x(t),x(\gamma (t))),$ & $t\neq t_{k}$ \\ 
$x(t_{k})-x\left(t_{k}^{-}\right) =J_{k}(x(t_{k}^{-})),$ & $t=t_{k},$ \\ 
$x(\tau )=x_{0},$ & 
\end{tabular}
\label{SISTEMA_GENERICO_IDEPCAG}
\end{equation}
where  $f\in C([\tau ,\infty )\times \mathbb{R}^{n}\times \mathbb{R}^{n},\mathbb{R}^{n}),$  $J_{k}\in C(\left\{t_{k}\right\} ,\mathbb{R}^{n})$ and $\left(\tau ,x_{0}\right) \in \mathbb{R}\times \mathbb{R}^{n}$.\\
Let the following hypothesis hold:
\begin{itemize}
\item[(H1)] Let $\eta_{1},\eta_{2}:\mathbb{R}\rightarrow [0,\infty )$ locally integrable functions 
and $\lambda_{k}\in \mathbb{R}^+$, $\forall k\in \mathbb{Z}$; such that
\begin{eqnarray*}
\left\Vert f(t,x_{1},y_{1})-f(t,x_{2},y_{2})\right\Vert &\leq &\eta_{1}(t) \left\Vert x_{1}-x_{2}\right\Vert 
+\eta_{2}(t) \left\Vert y_{1}-y_{2}\right\Vert ,\\
\left\Vert J_{k}(x_{1}(t_{k}^{-}))-J_{k}(x_{2}(t_{k}^{-}))\right\Vert 
&\leq &\lambda_{k}\left\Vert x_{1}\left( t_{k}^{-}\right) -x_{2}\left(t_{k}^{-}\right) \right\Vert .
\end{eqnarray*}
where $\Vert \cdot \Vert$ is some matricial norm.
\item[(H2)]
\begin{equation*}
\overline{\nu }=\sup_{k\in \mathbb{Z}}\left(\int_{t_{k}}^{t_{k+1}}\left(\eta_{1}(s)+\eta_{2}(s)\right)ds\right) <1.
\end{equation*}
\end{itemize}

In the following, we mention some useful results:
an integral equation associated with \eqref{sistema_idepcag_general_abstract} and two Gronwall-Bellman type inequalities necessary to prove the uniqueness and stability of 
solutions.
\subsection{An Integral equation associated to \eqref{SISTEMA_GENERICO_IDEPCAG}}
\begin{theorem}{(\cite{CTP2019}, Lemma 4.2)}
\label{TEOREMA_ECUACION_INTEGRAL}a function $x(t)=x(t,\tau ,x_0)$,  $\tau\in\mathbb{R}^+$ is a solution of \eqref{SISTEMA_GENERICO_IDEPCAG} on $\mathbb{R}^{+}$ if and only if satisfies: 
\begin{equation*}
x(t)=x_{0}+\int_{\tau}^{t}f(s,x(s),x(\gamma(s)) ds+\sum_{\tau\leq t_{k}<t}J_{k}\left(x\left(t_{k}^{-}\right) \right),
\end{equation*}
where
\begin{eqnarray*}
\int_{\tau}^{t}f(s,x(s),x\left(\gamma(t)\right)) ds &=&\int_{\tau}^{t_{1}}f(s,x(s),x\left(\zeta_{0}) \right) ds
+\sum_{j=1}^{k(t)-1}\int_{t_{j}}^{t_{j+1}}f(s,x(s),x(\zeta_{j}))ds \\
&&+\int_{t_{k(t)}}^{t}f\left(s,x(s),x\left(\zeta_{k(t)}\right)\right) ds,
\end{eqnarray*}
\end{theorem}

\subsection{\bf{First IDEPCAG Gronwall-Bellman type inequality}}
\begin{lemma}{(\cite{Torres1},\cite{CTP2019} Lemma 4.3)}
\label{LEMA_GRONWALL_PINTO_2}Let $I$ an interval and $u,\eta_{1},\eta_{2}:I\rightarrow [0,\infty )$ such that $u$ is continuous (with possible exception at $\{t_k\}_{k\in\mathbb{N}}$), $\eta_{1},\eta_{2}$ are continuous and locally integrable functions, $\eta:\left\{ t_{k}\right\} \rightarrow [0,\infty )$ and $\gamma (t) $ a piecewise constant argument of generalized type such that $\gamma (t)=\zeta_{k}$,  
$\forall t\in I_{k}=[t_{k},t_{k+1})$ with $t_{k}\leq \zeta_{k}\leq t_{k+1}$ $\forall k\in \mathbb{N}.$ 
Assume that $\forall t\geq \tau $  
\begin{equation*}
u(t)\leq u(\tau) +\int_{\tau}^{t}\left(\eta_{1}(s)u(s)+\eta_{2}(s)u(\gamma(s))\right) ds+\sum_{\tau \leq t_{k}<t}\eta(t_{k})u(t_{k}^{-})
\end{equation*}
and
\begin{equation}
\widehat{\vartheta }_{k}=\int_{t_{k}}^{\zeta_{k}}\left(\eta_{1}(s)+\eta_{2}(s)\right) ds\leq \widehat{\vartheta }
:=\sup_{k\in \mathbb{N}}\widehat{\vartheta }_{k}<1.  \label{CONDICION_INVERTIBILIDAD_GRONWALL_2}
\end{equation}
hold. Then, for $t\geq \tau $, we have
\begin{align}
u(t)\leq &\left(\prod_{\tau \leq t_{k}<t}\left(1+\eta(t_{k})\right) \right) 
\exp \left(\int_{\tau}^{t}\left(\eta_{1}(s)+\frac{\eta_{2}(s)}{1-\widehat{\vartheta}}\right) ds\right) u(\tau),\\
u(\zeta_k)\leq &(1-\vartheta)u(t_k)\\
u(\gamma(t))\leq &(1-\vartheta)^{-1}\left(\prod_{\tau \leq t_{k}<t}\left(1+\eta_3 (t_{j})\right) \right) 
\exp \left(\int_{\tau}^{t}\left(\eta_{1}(s)+\frac{\eta_{2}(s)}{1-\widehat{\vartheta }}\right) ds\right) u(\tau).
\label{DESIGUALDAD_GRONWALL_2}
\end{align}
\end{lemma}

\subsection{\bf{Second IDEPCAG Gronwall-Bellman type inequality}}
\begin{lemma}{(\cite{Torres1},\cite{kuo2022})}
\label{LEMA_GRONWALL_PINTO_1}Let $I$ an interval and $u,\eta_{1},\eta_{2}:I\rightarrow [0,\infty )$ such that $u$ is continuous (with possible exception at $\{t_k\}_{k\in N}$), $\eta_{1},\eta_{2}$ are continuous and locally integrable functions, $\eta:\left\{ t_{k}\right\} \rightarrow [0,\infty )$ and $\gamma(t)$ a piecewise constant argument of generalized type such that $\gamma (t)=\zeta_{k}$,  
$\forall t\in I_{k}=[t_{k},t_{k+1})$ with $t_{k}\leq \zeta_{k}\leq t_{k+1}$ $\forall k\in \mathbb{N}.$ 
Assume that $\forall t\geq \tau $  
\begin{equation}
u(t)\leq u(\tau) +\int_{\tau}^{t}(\eta_{1}(s)u(s)+\eta_{2}(s)u(\gamma(s))) ds
+\sum_{\tau \leq t_{k}<t}\eta(t_{k})u(t_{k}^{-})  \label{gronwall1}
\end{equation}
and
\begin{equation}
\varrho_{k}=\int_{t_{k}}^{\zeta_{k}}\left(\eta
_{2}(s)e^{\int_{s}^{\zeta_{k}}\eta_{1}(r)dr}\right) ds\leq \varrho
:=\sup_{k\in \mathbb{N}}\varrho_{k}<1.  \label{CONDICION_INVERTIBILIDAD_GRONWALL_1}
\end{equation}
Then, for $t\geq \tau $, we have
\begin{eqnarray}
u(t) &\leq &\left(\prod_{\tau \leq t_{k}<t}\left(1+\eta(t_{k})\right) \right) \notag\\
&&\cdot \exp \left(\frac{1}{1-\vartheta }\sum_{j=k(\tau)+1}^{k(t)}\int_{t_{j-1}}^{t_{j}}\right. 
\eta_{2}(s)\exp \left(\int_{t_{j-1}}^{\zeta_{j-1}}\eta_{1}(r)dr\right) ds  \label{DESIGUALDAD_GRONWALL_1} \\
&&+\frac{1}{1-\vartheta }\int_{t_{k(t)}}^{t}\eta_{2}(s)\exp\left. \left(\int_{t_{k(t)}}^{\zeta_{k(t)}}
\eta_{1}(r)dr\right) ds+\int_{\tau}^{t}\eta_{1}(s)ds\right) u(\tau).  \notag
\end{eqnarray}
\end{lemma}

\subsection{Existence and uniqueness for \eqref{SISTEMA_GENERICO_IDEPCAG}}
\begin{theorem}{(\bf{Uniqueness}) (\cite{CTP2019}, Theorem 4.5)}
Consider the I.V.P for \eqref{sistema_idepcag_general_abstract} with $y(t,\tau,y(\tau))$. Let (H1)-(H2) and \ref{LEMA_GRONWALL_PINTO_2} hold. Then, there exists a unique solution $y$ for \eqref{sistema_idepcag_general_abstract} on $[\tau,\infty)$. Moreover, every solution is stable.
\end{theorem}

\begin{lemma}{(\bf{Existence of solutions in $[\tau,t_{k})$}) (\cite{CTP2019}, Lemma 4.6)}
Consider the I.V.P for \eqref{sistema_idepcag_general_abstract} with $y(t,\tau,y(\tau))$. Let (H1)-(H2) and \ref{LEMA_GRONWALL_PINTO_2} hold. Then, for each $y_0\in\mathbb{R}^n$ and $\zeta_k\in [t_{k-1},t_{k})$ there exists a solution $y(t)=y(t,\tau,y(\tau))$ of \eqref{sistema_idepcag_general_abstract} on $[\tau,t_r)$ such that $y(\tau)=y_0$.
\end{lemma}

\begin{theorem}{(\bf{Existence of solutions in $[\tau,\infty$}) (\cite{CTP2019}, Theorem 4.7)}
Let (H1)-(H2) and \ref{LEMA_GRONWALL_PINTO_2} hold. Then, for each $(\tau,y_0\in\mathbb{R}_0^+\times\mathbb{R}^n$, there exists $y(t)=y(t,\tau,y_0)$ for $t\geq \tau$, a unique solution for \eqref{sistema_idepcag_general_abstract} such that $y(\tau)=y_0$.    
\end{theorem}

\section{Variation of parameters formula for IDEPCAG}
In this section, we will construct a variation of parameters formula for the IDEPCAG system
\begin{equation}
\begin{tabular}{ll}
$y^{\prime }(t)=A(t)y(t)+B(t)y(\gamma (t))+F(t),$ & $t\neq t_{k}$ \\ 
$\Delta y|_{t=t_{k}}=C_k y(t_{k}^{-})+D_k,$ & $t=t_{k}$
\end{tabular}
\label{SISTEMA_IDEPCAG_ORIGINAL_}
\end{equation}
where $y\in\mathbb{R}^{n\times 1},t\in\mathbb{R},F(t)\in\mathbb{R}^{n\times 1}$ is a real valued continuous matrix,  $A(t),B(t)\in \mathbb{R}^{n\times n}$ are real valued continuous locally integrable matrices, $C_k,D_k\in\mathbb{R}^{n\times n}$,  $(I+C_k)$ invertible $\forall k\in\mathbb{Z},$ where $I_{n\times n}=I$ is the identity matrix and $\gamma (t)$ is a generalized piecewise constant argument. This time, we will consider the advanced and the delayed intervals in our approach.\\

\noindent First, we will find the \textit{fundamental matrix} for the linear IDEPCAG
\begin{equation}
\begin{tabular}{ll}
$w^{\prime }(t)=A(t)w(t)+B(t)w(\gamma (t)),$ & $t\neq t_{k}$ \\ 
$\Delta w|_{t=t_{k}}=C_k w(t_{k}^{-}),$ & $t=t_{k}.$
\end{tabular}
\label{SISTEMA_IDEPCAG_ORIGINAL_HOMOGEGEO}
\end{equation}
Then, we will give the variation of parameters formula for \eqref{SISTEMA_IDEPCAG_ORIGINAL_}.\\

Let $\Phi(t,s), \,t,s\in\mathbb{R},$ with $\Phi(t,t)=I$ the transition (Cauchy) matrix of the ordinary system
\begin{equation}
\begin{tabular}{ll}
$x^{\prime }(t)=A(t)x(t),$ & $t\in I_{k}=[t_{k},t_{k+1}).$
\end{tabular}
\   \label{sistema_homogeneo_ordinario}
\end{equation}

We will assume the following hypothesis:
\begin{enumerate}
    \item[(H3)] Let
\begin{eqnarray*}
\rho_{k^{+}}(A) &=&exp\left(\int_{t_{k}}^{\zeta_{k}}\left\Vert A(u)\right\Vert du\right),\qquad \rho_{k^{-}}(A)
=exp\left(\int_{\zeta_{k}}^{t_{k+1}}\left\Vert A(u)\right\Vert du\right), \\
\rho_{k}(A) &=&\rho_{k^{+}}(A)\cdot \rho_{k^{-}}(A),\qquad  \nu_k^{\pm}(B)=\rho_k^{\pm }(A)\ln \rho_k^{\pm }(B).
\end{eqnarray*}
and assume that $$\rho(A)=\displaystyle{\sup_{k\in \mathbb{Z}}\rho_{k}(A)<\infty},\qquad \nu^{\pm}(B)=\sup_{k\in \mathbb{Z}}\nu_{k}^{\pm}(B)<\infty.$$
Consider the following matrices
\begin{equation}
J(t,\tau)=I+\int_{\tau}^{t}\Phi (\tau,s)B(s)ds,\qquad   
E(t,\tau)=\Phi(t,\tau)J(t,\tau), \label{MATRIZ J}
\end{equation}
where 
\begin{equation}
\nu_k^{\pm}(B)<\nu^{\pm}(B)<1. \label{Cond_invertibilidad}
\end{equation}
\end{enumerate}

\begin{remark}
It is important to notice the following facts: 
\begin{enumerate}
\item[a)] As a consequence of (H3), $J(t_k,\zeta_k)$ and $J(t_{k+1},\zeta_k)$ are invertible $\forall
k\in \mathbb{Z},$ 
and 
\begin{eqnarray*}
\left\Vert J^{-1}(t_k,\zeta_{k})\right\Vert \leq \sum\limits_{k=0}^{\infty } \left[ \nu^{+} (B)\right] ^{k}  =\frac{1}{1-\nu^{+} (B)},\quad 
\left\Vert J(t_k,\zeta_{k})\right\Vert \leq 1+\nu^+ (B),  \label{COTA_J}\\
\left\Vert J^{-1}(t_{k+1},\zeta_{k})\right\Vert \leq \sum\limits_{k=0}^{\infty } \left[ \nu^{-} (B)\right] ^{k}  =\frac{1}{1-\nu^{-} (B)},\quad 
\left\Vert J(t_{k+1},\zeta_{k})\right\Vert \leq 1+\nu^- (B).  \label{COTA_J2}
\end{eqnarray*}
Additionally, setting $t_0\coloneqq \tau$ we will assume that $J^{-1}(\tau,\gamma(\tau))$ exists.
\item [b)] Also, due to (H3) and the Gronwall inequality, we have 
\begin{equation*}
    |\Phi(t)|\leq \rho(A),
\end{equation*}
(See \cite{P2011}).
\end{enumerate}
\end{remark}

\subsection{The fundamental matrix of the linear homogeneous IDEPCAG}~\par
\vskip1mm
We adopt the following convention:
\begin{equation*}
^{\leftarrow }\prod_{k=j}^{j+p}T_{k}=T_{j+p}\cdot T_{j+p-1}\cdot \ldots\cdot T_{j}.
\end{equation*}
Also, we will assume $\gamma(\tau)\coloneqq\tau$ if $\gamma(\tau)<\tau$, where $k(\tau)$ is the only $k\in\mathbb{Z}$ such that $\tau \in I_{k(\tau)}=[t_{k(\tau)},t_{k(\tau)+1}).$ 
We will adopt the following notation:
$$\prod_{j=r+1}^r A_j=1,\qquad \sum_{j=r+1}^r A_j=0.$$
Let the system
\begin{equation}
\begin{tabular}{ll}
$w^{\prime }(t)=A(t)w(t)+B(t)w(\gamma (t)),$ & $t\neq t_{k}$ \\ 
$w(t_{k})=\left( I+C_k\right) w(t_{k}^{-}),$ & $t=t_{k}$\\
$w_0=w(\tau).$
\end{tabular}
\label{sistema_w}
\end{equation}
We will construct the fundamental matrix for system \eqref{sistema_w}.\\
Let $t,\tau \in I_k=[t_{k},t_{k+1})$ for some $k\in \mathbb{Z}.$
In this interval, we are in the presence of the ordinary system
$$w^{\prime }(t)=A(t)w(t)+B(t)w(\zeta_k).$$
So, the unique solution can be written as
\begin{equation}
w(t)=\Phi(t,\tau)w(\tau )+\int_{\tau}^{t}\Phi(t,s)B(s)w(\zeta_{k})ds.
\label{variacion_parametros_general}
\end{equation}
Keeping in mind $I_k^{+}=[t_k,\zeta_k]$, evaluating the last expression at  $t=\zeta_{k}$ we have
\begin{equation}
w(\zeta_{k})=\Phi (\zeta_{k},\tau )w(\tau)+\int_{\tau}^{\zeta_{k}}\Phi(\zeta_{k},s)B(s)w(\zeta_{k})ds,
\label{variacion_parametros_gamma_general}
\end{equation}
Hence, we get
\begin{eqnarray*}
\left( I+\int_{\zeta_{k}}^{\tau }\Phi (\zeta_{k},s)B(s)ds\right) w(\zeta_{k}) &=&\Phi (\zeta_{k},\tau )w(\tau).
\end{eqnarray*}
I.e
\begin{equation}
w(\zeta_{k})=J^{-1}(\tau ,\zeta_{k})\Phi (\zeta_{k},\tau )w(\tau).
\label{inicial_gamma_general}
\end{equation}
Then, by the definition of $E(t,\tau)=\Phi(t,\tau)J(t,\tau)$, we have
\begin{equation}
w(\zeta_{k})=E^{-1}(\tau ,\zeta_{k})w(\tau).
\label{condicion_inicial_con_gama_y_E}
\end{equation}

Now, from \eqref{variacion_parametros_general} working on $I_k^{-}=[\zeta_k,t_{k+1})$, considering $\tau =\zeta_{k}$, we have
\begin{eqnarray*}
w(t) &=&\Phi (t,\zeta_{k})w(\zeta_{k})+\int_{\zeta_{k}}^{t}\Phi(t,s)B(s)w(\zeta_{k})ds \\
&=&\Phi (t,\zeta_{k})\left( I+\int_{\zeta_{k}}^{t}\Phi (\zeta_{k},s)B(s)ds\right) w(\zeta_{k}).
\end{eqnarray*}
I.e.,
\begin{equation}
w(t)=E(t,\zeta_{k})w(\zeta_{k}).
\label{lineal_homogenea_parte_retardo_sin_cond_inicial}
\end{equation}

So, by \eqref{condicion_inicial_con_gama_y_E}, we can rewrite  \eqref{lineal_homogenea_parte_retardo_sin_cond_inicial} as
\begin{equation}
w(t)=E(t,\zeta_{k})E^{-1}(\tau ,\zeta_{k})w(\tau).
\label{condicion_inicial_sin_w}
\end{equation}

Then, setting  
\begin{equation}
W(t,s)=E(t,\gamma(s))E^{-1}(s ,\gamma(s)),\qquad \text{if }t,s
\in I_k=[t_{k},t_{k+1}),  \label{MATRIZ_w}
\end{equation}
we have the solution for \eqref{sistema_w} for $t\in I_{k}$
\begin{equation}
w(t)=W(t,\tau)w(\tau). \label{DEFINICION_GENERAL_MATRIZ_W}
\end{equation}

Next, if we consider $\tau=t_k $ and assuming left side continuity of \eqref{DEFINICION_GENERAL_MATRIZ_W} at $t=t_{k+1}$, we have
\begin{equation*}
w(t_{k+1}^{-})=W(t_{k+1},t_k )w(t_k)
\end{equation*}
Then, applying the impulsive condition to the last equation, we get
\begin{eqnarray*}
w(t_{k+1}) &=&\left( I+C_{k+1}\right) w(t_{k+1}^{-}) \\
&=&\left( I+C_{k+1}\right) W(t_{k+1},t_k)w(t_k).
\end{eqnarray*}

This expression corresponds to a finite-difference equation. Then, by solving it, we get
\begin{equation}
w(t_{k(t)})=\left(^{\leftarrow}\prod_{j=k(\tau)+1}^{k(t)-1}\left( I+C_{j+1}\right) W(t_{j+1},t_{j})\right) w\left( t_{k(\tau)+1}\right).
\label{SOLUCION_DISCRETA_SISTEMA_Z}
\end{equation}

Finally, by \eqref{DEFINICION_GENERAL_MATRIZ_W} and the impulsive condition, we have
$$w(t_{k(\tau)+1})=(I+C_{k(\tau)+1})W(t_{k(\tau)+1},\tau)w(\tau).$$
Hence, considering $\tau=t_k$ in \eqref{DEFINICION_GENERAL_MATRIZ_W} and applying \eqref{SOLUCION_DISCRETA_SISTEMA_Z} 
we get
\begin{eqnarray}
w(t) &=&W(t,t_{k(t)})\left(^{\leftarrow}
\prod_{j=k(\tau)+1}^{k(t)-1}\left(I+C_{j+1}\right) W(t_{j+1},t_{j})\right)
\left( I+C_{k(\tau)+1}\right) W(t_{k(\tau)+1},\tau)w(\tau)
 \notag\\
&=&W(t,\tau)w(\tau),\qquad \text{for }t\in I_{k(t)}\text{ and }\tau \in I_{k(\tau) }.  \label{SOLUCION_FINAL_SISTEMA_IDEPCAG_LINEAL}
\end{eqnarray}
The last equation is the solution of \eqref{sistema_w} on $[\tau,t).$\\
We call to the expression 
\begin{equation}
W(t,\tau)=W(t,t_{k(t)})\left(^{\leftarrow}\prod_{j=k(\tau)+1}^{k(t)-1}\left( I+C_{j+1}\right) W(t_{j+1},t_{j})\right) 
\left( I+C_{k\left(\tau \right) +1}\right) W(t_{k(\tau)+1},\tau ),
\label{MATRIZ_FUNDAMENTAL_IDEPCAG}
\end{equation}
the fundamental matrix for \eqref{sistema_w}
for $t\in I_{k(t)}$ and $\tau \in I_{k(\tau) }.$

\begin{remark}\label{decomposition_interval_W}
We use the decomposition of $I_k=I_k^+\cup I_k^-$ to define $W$. In fact, we can rewrite \eqref{MATRIZ_FUNDAMENTAL_IDEPCAG} in terms of the advanced and delayed parts using \eqref{MATRIZ_w}:
\begin{eqnarray*}
W(t,\tau)&=&E(t,\zeta_{k(t)})E^{-1}(t_{k(t)},\zeta_{k(t)})\left(^{\leftarrow}\prod_{j=k(\tau)+1}^{k(t)-1}\left(I+C_{j+1}\right) E(t_{j+1},\zeta_{j})E^{-1}(t_{j},\zeta_{j})\right)\\ 
&&\cdot \left( I+C_{k\left(\tau \right) +1}\right) E(t_{k(\tau)+1},\gamma(\tau))E^{-1}(\tau,\gamma(\tau)),\qquad \zeta_j=\gamma(t_j).
\label{MATRIZ_FUNDAMENTAL_IDEPCAG_avance_retardo}
\end{eqnarray*}
for $t\in I_{k(t)}$ and $\tau \in I_{k(\tau) }.$
\end{remark}
\bigskip
\bigskip
\bigskip
\begin{remark}
\item[a)]   Considering $B(t)=0$, we recover the classical fundamental matrix of the impulsive linear differential equation (see \cite{Samoilenko}).
\item[b)] If $C_k=0, \forall k\in\mathbb{Z}$, we recover the DEPCAG case studied by M. Pinto in \cite{P2011}. 
\item[c)] If we consider $\gamma(t)=p\left[\dfrac{t+l}{p}\right],$ with $p<l$, we recover the IDEPCA case studied by K-S. Chiu in \cite{Kuo2023}. 
\end{remark}
\bigskip
\subsection{The variation of parameter formula for IDEPCAG}~\par
\vskip1mm
Let the IDEPCAG
\begin{equation}
\begin{tabular}{ll}
$y^{\prime }(t)=A(t)y(t)+B(t)y(\gamma (t))+F(t),$ & $t\neq t_{k}$ \\ 
$y(t_{k})=(I+{C}_{k})y(t_{k}^{-})+D_k,$ & $t=t_{k}$ \\ 
$y_0=y(\tau). $ & 
\end{tabular}
\label{sistema_idepcag_lineal_no_homogeneo}
\end{equation}
If $\tau,t\in I_{k}=[t_k,t_{k+1})$, then the unique solution of \eqref{sistema_idepcag_lineal_no_homogeneo} is 
\begin{equation*}
y(t)=\Phi(t,\tau)y(\tau)+\int_{\tau}^{t}\Phi(t,s)B(s)y(\zeta_{k})ds+\int_{\tau}^{t}\Phi(t,s)f(s)ds.
\end{equation*}
Then, if $\tau =\zeta_{k}$, we have
\begin{eqnarray*}
y(t) &=&\Phi (t,\zeta_{k})y(\zeta_{k})+\int_{\zeta_{k}}^{t}\Phi(t,s)B(s)y(\zeta_{k})ds+\int_{\zeta_{k}}^{t}\Phi(t,s)f(s)ds \\
&=&\Phi (t,\zeta_{k})\left( I+\int_{\zeta_{k}}^{t}\Phi (\zeta_{k},s)B(s)ds\right) y(\zeta_{k})+\int_{\zeta_{k}}^{t}\Phi(t,s)f(s)ds \\
&=&\Phi (t,\zeta_{k})J\left( t,\zeta_{k}\right) y(\zeta_{k})+\int_{\zeta_{k}}^{t}\Phi(t,s)f(s)ds,
\end{eqnarray*}
I.e
\begin{equation}
y(t)=E(t,\zeta_{k})y(\zeta_{k})+\int_{\zeta_{k}}^{t}\Phi(t,s)f(s)ds.
\label{variacion_parametros_con_E_sin_cond_inicial}
\end{equation}

Now, if we consider $t=\tau $ in \eqref{variacion_parametros_con_E_sin_cond_inicial}
we have
\begin{equation*}
y(\tau)=E(\tau ,\zeta_{k})y(\zeta_{k})+\int_{\zeta_{k}}^{\tau }\Phi(\tau ,s)f(s)ds,
\end{equation*}
and, by $(H3)$, we get the following estimation for $y(\zeta_k)$
\begin{equation}
y(\zeta_{k})=E^{-1}(\tau ,\zeta_{k})\left( y(\tau)+\int_{\tau}^{\zeta_{k}}\Phi(\tau,s)f(s)ds\right) .
\label{condicion_inicial_gamma_caso_lineal_no_homog}
\end{equation}

Then, applying \eqref{condicion_inicial_gamma_caso_lineal_no_homog} in 
\eqref{variacion_parametros_con_E_sin_cond_inicial} we obtain
\begin{equation*}
y(t)=E(t,\zeta_{k})E^{-1}(\tau ,\zeta_{k})\left( y(\tau)+\int_{\tau}^{\zeta_{k}}\Phi(\tau,s)f(s)ds\right) +\int_{\zeta_{k}}^{t}\Phi(t,s)f(s)ds,
\end{equation*}

i.e.,%
\begin{equation}
y(t) =W(t,\tau)y(\tau)+\int_{\tau}^{\zeta_{k}}W(t,\tau)\Phi(\tau ,s)f(s)ds+\int_{\zeta_{k}}^{t}\Phi(t,s)f(s)ds,\quad \tau ,t\in I_{k}.
\label{solucion_sin_cond_inicial_en_terminos_de_W}
\end{equation}

Next, taking the left-side limit to the last expression, we have
\begin{equation}
y(t_{k+1}^{-})=W(t_{k+1},\tau )\left( y(\tau)+\int_{\tau}^{\zeta_{k}}\Phi(\tau,s)f(s)ds\right) +\int_{\zeta_{k}}^{t_{k+1}}\Phi (t_{k+1},s)f(s)ds.
\label{limite_lateral_caso_idepcag_lineal_no_homog}
\end{equation}

Applying the impulsive condition, we get
\begin{equation}
y(t_{k+1}) =\left( I+C_{k+1}\right) y(t_{k+1}^{-})+D_{k+1},\notag
\end{equation}
or
\begin{eqnarray*}
y(t_{k+1}) &=&\left( I+C_{k+1}\right) W(t_{k+1},\tau )\left(y(\tau)+\int_{\tau}^{\zeta_{k}}\Phi(\tau,s)f(s)ds\right) \\
&&+\int_{\zeta_{k}}^{t_{k+1}}\left( I+C_{k+1}\right) \Phi(t_{k+1},s)f(s)ds+D_{k+1},
\end{eqnarray*}
Therefore, considering  $\tau =t_{k}$ in the last expression we have
\begin{eqnarray*}
y(t_{k+1}) &=&\left( I+C_{k+1}\right) W(t_{k+1},t_k)\left(y(t_k)+\int_{t_k}^{\zeta_{k}}\Phi(t_k,s)f(s)ds\right) \\
&&+\int_{\zeta_{k}}^{t_{k+1}}\left( I+C_{k+1}\right) \Phi(t_{k+1},s)f(s)ds+D_{k+1},
\end{eqnarray*}
or
\begin{eqnarray*}
y(t_{k+1}) &=&W_k\left(y(t_k)+\alpha_k^+\right)+\alpha_k^- +\beta_k,
\end{eqnarray*}
which corresponds to a non-homogeneous linear difference equation,  
where 
\begin{align*}
&W_k=(I+C_{k+1})W(t_{k+1},t_k),\\
&\alpha_k^+=\int_{t_k}^{\zeta_{k}}\Phi(t_k,s)f(s)ds,\\
&\alpha_k^-=\int_{\zeta_{k}}^{t_{k+1}}\left( I+C_{k+1}\right) \Phi(t_{k+1},s)f(s)ds,\\
&\beta_k=D_{k+1}.    
\end{align*}
Recalling that 
\begin{equation*}
W(t_{k(t)},\tau)=\left(^{\leftarrow}\prod_{j=k(\tau)+1}^{k(t)-1}\left( I+C_{j+1}\right) W(t_{j+1},t_{j})\right) 
\left( I+C_{k\left(\tau \right) +1}\right) W(t_{k(\tau)+1},\tau ),
\end{equation*}
we get the discrete solution of \eqref{sistema_idepcag_lineal_no_homogeneo}:
\begin{eqnarray*}
y(t_{k(t)}) &=&\left(^{\leftarrow}\prod_{j=k(\tau) +1}^{k(t)-1}\left( I+C_{j+1}\right) W(t_{j+1},t_{j})\right)
(I+C_{k(\tau)+1})W(t_{k(\tau)+1},\tau)y(\tau)  \notag \\
&&+\int_{\tau}^{\zeta_{k(\tau)}}W\left(t_{k(t)},\tau \right) \Phi(\tau,s)f(s)ds  \notag \\
&&+\sum_{r=k(\tau)+1}^{k(t)-1}\left(^{\leftarrow}\prod_{j=r}^{k(t)-1}\left(I+C_{j+1}\right)W(t_{j+1},t_{j})\right)\int_{t_{r}}^{\zeta_{r}}\Phi (t_{r},s)f(s)ds  \notag \\
&&+\sum_{r=k(\tau)}^{k(t)-1}\left(^{\leftarrow}\prod_{j=r+1}^{k(t)-1}(I+C_{j+1}) W(t_{j+1},t_{j})\right) \int_{\zeta_{r}}^{t_{r+1}}(I+C_{r+1}) \Phi(t_{r+1},s)f(s)ds
\notag \\
&&+\sum_{r=k(\tau)}^{k(t)-1}\left(^{\leftarrow}\prod_{j=r+1}^{k(t)-1}(I+C_{j+1}) W(t_{j+1},t_{j})\right) D_{r+1}.  \notag
\end{eqnarray*}
or, written in terms of \eqref{MATRIZ_FUNDAMENTAL_IDEPCAG},
\begin{eqnarray}
y(t_{k(t)}) &=&W(t_{k(t)},\tau)y(\tau)+\int_{\tau}^{\zeta_{k(\tau) }}W\left(t_{k(t)},\tau \right) \Phi(\tau,s)f(s)ds  \label{solucion_discreta_idepcag_lineal_no_homog} \\
&&+\sum_{r=k(\tau)+1}^{k(t)-1}\int_{t_{r}}^{\zeta_{r}}W(t_{k(t)},t_{r})\Phi(t_{r},s)f(s)ds  \notag \\
&&+\sum_{r=k(\tau)}^{k(t)-1}\int_{\zeta_{r}}^{t_{r+1}}W(t_{k(t)},t_{r+1})(I+C_{r+1}) \Phi(t_{r+1},s)f(s)ds
\notag \\
&&+\sum_{r=k(\tau)}^{k(t)-1}W(t_{k(t)},t_{r+1})D_{r+1}.  \notag
\end{eqnarray}
Now, considering $\tau=t_k$ in \eqref{solucion_sin_cond_inicial_en_terminos_de_W}  we have 
\begin{eqnarray*}
    y(t)&=&W(t,t_{k(t)})y(t_{k(t)})\\
&&+\int_{t_{k(t)}}^{\zeta_{k(t)}}W(t,t_{k(t)})\Phi(t_{k(t)},s)f(s)ds+\int_{\zeta_{k(t)}}^t\Phi(t,s)f(s)ds.
\end{eqnarray*}
Finally, replacing $y(t_{k(t)})$ by \eqref{solucion_discreta_idepcag_lineal_no_homog} and rewriting in terms of 
\eqref{MATRIZ_FUNDAMENTAL_IDEPCAG}, we get the variation of parameters formula for IDEPCAG \eqref{sistema_idepcag_lineal_no_homogeneo}:
\begin{eqnarray}
y(t)&=&W(t,\tau)y(\tau)
\label{FORMULA_VARIACION_PARAMETROS_IDEPCAG_GENERAL_FINAL} \\
&&+\int_{\tau}^{\zeta_{k(\tau) }}W(t,\tau) \Phi(\tau ,s)f(s)ds
+\sum_{r=k(\tau)+1}^{k(t)}\int_{t_{r}}^{\zeta_{r}}W(t,t_{r})\Phi(t_{r},s)f(s)ds  \notag \\
&&+\sum_{r=k(\tau) }^{k(t)-1}\int_{\zeta_{r}}^{t_{r+1}}W(t,t_{r+1})\left( I+C_{r+1}\right) \Phi(t_{r+1},s)f(s)ds  \notag \\
&&+\int_{\zeta_{k(t)}}^{t}\Phi(t,s)f(s)ds+\sum_{r=k(\tau)+1}^{k(t)}W(t,t_{r})D_{r}, \quad \text{for }t\in [\tau ,t_{k(t)+1}),  \notag
\end{eqnarray}
where $W$ is the fundamental matrix of \eqref{sistema_w}.\\

\bigskip
\subsubsection{Green type matrix for IDEPCAG}\hfill\\
If we define the following Green matrix type for IDEPCAG: 
\begin{equation}
\widetilde{W}(t,s)=\left\{ 
\begin{array}{l}
W^{+}(t,s),\text{ \ if }t_{r}\leq s\leq \gamma (s) \\ 
\\ 
W^{-}(t,s),\text{ \ if }\gamma (s) <s\leq t_{r+1},
\end{array}
\right.  \label{GREEN_FORMULA:VARIACION_PARAMETROS_IDEPCAG}
\end{equation}
where
\begin{equation}
W^{+}(t,s)= 
\begin{array}{c}
W(t,t_{r})\Phi(t_{r},s)
\end{array}
\begin{array}{c}
\text{if \ }t_{r}\leq s\leq \gamma (s) ,\text{ }s<t,
\end{array}
\label{GREEN_VARIACION_PARAMETROS_MAS}
\end{equation}
and
\begin{equation}
W^{-}(t,s)=\left\{ 
\begin{array}{c}
W(t,t_{r+1})\left( I+C_{r+1}\right) \Phi (t_{r+1},s) \\ 
\Phi(t,s)
\end{array}
\begin{array}{l}
\text{if }\gamma (s) \leq s<t_{r+1},\text{ }t>s \\ 
\text{if }\gamma (t) <s\leq t<t_{r+1}.
\end{array}
\right.   \label{GREEN_VARIACION_PARAMETROS_MENOS}
\end{equation}
Hence, we can see that
\begin{eqnarray*}
\int_{\tau}^{t}W^{+}(t,s)f(s)ds&=&\int_{\tau}^{\zeta_{k(\tau) }}W(t,\tau) \Phi(\tau,s)f(s)ds\\
&&+\sum_{r=k(\tau)+1}^{k(t)}\int_{t_{r}}^{\zeta_{r}}W(t,t_{r})\Phi(t_{r},s)f(s)ds, 
\end{eqnarray*}
\begin{eqnarray*}
\int_{\tau}^{t}W^{-}(t,s)f(s)ds&=&\sum_{r=k(\tau)}^{k(t)-1}\int_{\zeta_{r}}^{t_{r+1}}W(t,t_{r+1})\left( I+C_{r+1}\right)
\Phi (t_{r+1},s)f(s)ds\\
&&+\int_{\zeta_{k(t)}}^{t}\Phi(t,s)f(s)ds.
\end{eqnarray*}
So, we have
\begin{equation*}
\widetilde{W}(t,s)=W^{+}(t,s) +W^{-}(t,s) .
\end{equation*}
In this way, \eqref{FORMULA_VARIACION_PARAMETROS_IDEPCAG_GENERAL_FINAL} can be expressed as
\begin{equation}
y(t)=W(t,\tau)y(\tau)+\int_{\tau}^{t}\widetilde{W}(t,s)
f(s) ds+\sum_{r=k(\tau)+1}^{k(t)}W(t,t_{r})D_{r}.
\label{FORMULA_VARIACION_PARAMETROS_IDEPCAG_GREEN}
\end{equation}
\subsection{Some special cases of \eqref{sistema_idepcag_lineal_no_homogeneo}}\hfill\\
In the following, we present some r cases for \eqref{sistema_idepcag_lineal_no_homogeneo}.\\

\begin{enumerate}
\item Let $\gamma^-(t)=t_k$ and $\gamma^+(t)=t_{k+1}$, for all $t\in I_k=[t_k,t_{k+1}).$ I.e., we are considering the completely delayed and advanced general piecewise constant arguments. Then, taking in account Remark \ref{decomposition_interval_W}, the solution of \eqref{sistema_idepcag_lineal_no_homogeneo} for both cases $y_-(t)$ and $y_+(t)$ respectively are:
\begin{eqnarray}
&&y_-(t)=W_-(t,\tau)y(\tau)
\label{FORMULA_VARIACION_PARAMETROS_IDEPCAG_GENERAL_FINAL_DELAYED}\\
&&+\sum_{r=k(\tau) }^{k(t)-1}\int_{t_{r}}^{t_{r+1}}W_-(t,t_{r+1})\left( I+C_{r+1}\right) \Phi(t_{r+1},s)f(s)ds\notag \\
&&+\int_{t_{k(t)}}^{t}\Phi(t,s)f(s)ds +\sum_{r=k(\tau)+1}^{k(t)}W_-(t,t_{r})D_{r},  \notag
\end{eqnarray}
where
\begin{eqnarray*}
W_-(t,\tau)&=&E(t,t_{k(t)})\left(^{\leftarrow}\prod_{j=k(\tau)+1}^{k(t)-1}\left(I+C_{j+1}\right) E(t_{j+1},t_{j})\right)\\ 
&&\cdot \left( I+C_{k\left(\tau \right)+1}\right) E(t_{k(\tau)+1},\tau),
\label{MATRIZ_FUNDAMENTAL_IDEPCAG_retardo}
\end{eqnarray*}
and
\begin{eqnarray}
&&y_+(t)=W_+(t,\tau)y(\tau)
\label{FORMULA_VARIACION_PARAMETROS_IDEPCAG_GENERAL_FINAL_advanced} \\
&&+\int_{\tau}^{t_{k(\tau)+1 }}W_+(t,\tau) \Phi(\tau ,s)f(s)ds
+\sum_{r=k(\tau)+1}^{k(t)}\int_{t_{r}}^{t_{r+1}}W_+(t,t_{r})\Phi(t_{r},s)f(s)ds  \notag \\
&&-\int_{t}^{t_{k(t)+1}}\Phi(t,s)f(s)ds +\sum_{r=k(\tau)+1}^{k(t)}W_+(t,t_{r})D_{r},\notag
\end{eqnarray}
where 
\begin{eqnarray*}
W_+(t,\tau)&=&E(t,t_{k(t)+1})E^{-1}(t_{k(t)},t_{k(t)+1})\left(^{\leftarrow}\prod_{j=k(\tau)+1}^{k(t)-1}\left(I+C_{j+1}\right) E^{-1}(t_{j},t_{j+1})\right)\\ 
&&\cdot \left( I+C_{k\left(\tau \right) +1}\right) E^{-1}(\tau,t_{k(\tau)+1}),
\label{MATRIZ_FUNDAMENTAL_IDEPCAG_avance}
\end{eqnarray*}
for $t\in I_{k(t)}$ and $\tau \in I_{k(\tau) },$
recalling that  $\gamma(\tau)\coloneqq\tau$ if 
$\gamma(\tau)<\tau$.\\

\item Let the IDEPCAG
\begin{equation}
\begin{tabular}{ll}
$w^{\prime }(t)=B(t)w(\gamma (t)),$ & $t\neq t_{k}$ \\ 
$w(t_{k})=(I+{C}_{k})w(t_{k}^{-}),$ & $t=t_{k}$ \\ 
$w_0=w(\tau). $ & 
\end{tabular}
\label{sistema_idepcag_lineal_homogeneo_B_cero}
\end{equation}
We see that $\Phi(t,s)=I,\,$$E(t,s)=J(t,s)$ and $J(t,s)=I+\int_s^t B(u)du,$ where $I$ is the identity matrix. Hence the fundamental matrix for \eqref{sistema_idepcag_lineal_homogeneo_B_cero} is given by
\begin{eqnarray*}
W(t,\tau)&=&J(t,\zeta_{k(t)})J^{-1}(t_{k(t)},\zeta_{k(t)})\left(^{\leftarrow}\prod_{j=k(\tau)+1}^{k(t)-1}\left(I+C_{j+1}\right) J(t_{j+1},\zeta_{j})J^{-1}(t_{j},\zeta_{j})\right)\\ 
&&\cdot \left( I+C_{k\left(\tau \right) +1}\right) J(t_{k(\tau)+1},\gamma(\tau))J^{-1}(\tau,\gamma(\tau)),\qquad \zeta_j=\gamma(t_j).
\label{MATRIZ_FUNDAMENTAL_IDEPCAG_B_cero}
\end{eqnarray*}
for $t\in I_{k(t)}$ and $\tau \in I_{k(\tau) }.$\\
This case is very important because it is used for the approximation of solutions of differential equations considering $\gamma(t)=\left[\frac{t}{h}\right]h,$ with $h>0$ fixed.\\

\item Let the IDEPCAG
\begin{equation}
\begin{tabular}{ll}
$w^{\prime }(t)=Aw(t)+Bw(\gamma (t)),$ & $t\neq t_{k}$ \\ 
$w(t_{k})=(I+C)w(t_{k}^{-}),$ & $t=t_{k}$ \\ 
$w_0=w(\tau), $ & 
\end{tabular}
\label{sistema_idepcag_lineal_homogeneo_constante}
\end{equation}
and
\begin{equation}
\begin{tabular}{ll}
$y^{\prime }(t)=Ay(t)+By(\gamma (t))+f(t),$ & $t\neq t_{k}$ \\ 
$y(t_{k})=(I+C)y(t_{k}^{-})+D_k,$ & $t=t_{k}$ \\ 
$y_0=y(\tau), $ & 
\end{tabular}
\label{sistema_idepcag_lineal_no_homogeneo_constante}
\end{equation}
where $A^{-1}$ exist. By $(H3)$, we know that $J(t,\tau)=I+\int_\tau^t e^{A(\tau-s)}B\,\, ds$ is invertible, for $\tau,t\in I_{k}=[t_k,t_{k+1})$. Moreover, following \cite{P2011}, we see that
\begin{eqnarray}
&J(t,\tau)&=I+\int_\tau^t e^{A(\tau-s)}B ds\notag\\
&&=I+e^{A\tau}\left(\int_\tau^t (-A)e^{-As} ds\right) (-A^{-1})B\notag\\
&&=I+A^{-1}\left(I-e^{A(\tau-t)}\right)B.
\end{eqnarray}
Then, as $E(t,\tau)=\Phi(t,\tau)J(t,\tau)$, we have
\begin{equation}
    E(t,\tau)=e^{A(t-\tau)}\left(I+A^{-1}\left(I-e^{-A(t-\tau)}\right)B\right).
\end{equation}
In light of the last calculations, we define
\begin{eqnarray*}
&&\widetilde{E}(t)=e^{At}\left(I+A^{-1}\left(I-e^{-At}\right)B\right)\\
&&\eta^+_k=\zeta_k-t_k,\qquad \eta^-_k=t_{k+1}-\zeta_k,\quad k\in\mathbb{Z},\\
&&\eta(t)=t-\gamma(t).
\end{eqnarray*}
Recalling that
\begin{equation}
\widehat{W}(t,s)=\widetilde{E}(t-\gamma(s))\widetilde{E}^{-1}(\eta(s)),\qquad \text{if }t,s
\in I_k=[t_{k},t_{k+1}), 
\end{equation}
the solution of \eqref{sistema_idepcag_lineal_homogeneo_constante} is
$$w(t)=\widehat{W}(t,\tau)w(\tau),$$
where 
\begin{eqnarray*}
\widehat{W}(t,\tau)&=&\widetilde{E}(\eta(t))\widetilde{E}^{-1}(-\eta^+_{k(t)})\left(^{\leftarrow}\prod_{j=k(\tau)+1}^{k(t)-1}(I+C) \widetilde{E}(\eta^-_{j})\widetilde{E}^{-1}(-\eta^+_{j})\right)\\ 
&&\cdot ( I+C) \widetilde{E}(\eta^-_{k(\tau)+1})\widetilde{E}^{-1}(\eta(\tau)),
\label{MATRIZ_FUNDAMENTAL_IDEPCAG_constante}
\end{eqnarray*}
is the fundamental matrix for \eqref{sistema_idepcag_lineal_homogeneo_constante} with $t\in I_{k(t)}$ and $\tau \in I_{k(\tau) }.$\\
The solution for \eqref{sistema_idepcag_lineal_no_homogeneo_constante} is given by
\begin{eqnarray}
y(t) &=&\widetilde{E}(\eta(t))\widetilde{E}^{-1}(-\eta^+_{k(t)})\left(^{\leftarrow}\prod_{j=k(\tau)+1}^{k(t)-1}(I+C) \widetilde{E}(\eta^-_{j})\widetilde{E}^{-1}(-\eta^+_{j})\right)\notag\\ 
&&\cdot ( I+C) \widetilde{E}(\eta^-_{k(\tau)+1})\widetilde{E}^{-1}(\eta(\tau))\left(y(\tau)+\int_{\tau}^{\zeta_{k(\tau) }} e^{A(\tau-s)}f(s)ds\right)  \label{FORMULA_VARIACION:PARAMETROS_IDEPCAG_GENERAL_sistema_constante} \\
&&+\widetilde{E}(\eta(t))\widetilde{E}^{-1}(-\eta^+_{k(t)})\notag\\
&&\cdot \left\{\sum_{r=k(\tau) +1}^{k(t)}\left(^{\leftarrow}\prod_{j=r}^{k(t)-1}(I+C)
\widetilde{E}(\eta^-_{j})\widetilde{E}^{-1}(-\eta^+_{j})\right) \int_{t_{r}}^{\zeta_{r}}e^{A(t_{r}-s)}f(s)ds\right. \notag \\
&&+\sum_{r=k(\tau)}^{k(t)-1}\left(^{\leftarrow}\prod_{j=r+1}^{k(t)}(I+C)
\widetilde{E}(\eta^-_{j})\widetilde{E}^{-1}(-\eta^+_{j})\right) \int_{\zeta_{r}}^{t_{r+1}}(I+C)e^{A(t_{r+1}-s)}f(s)ds  \notag \\
&&+\left.\sum_{r=k(\tau)}^{k(t)-1}\left(^{\leftarrow}\prod_{j=r+1}^{k(t)}(I+C)\widetilde{E}(\eta^-_{j})\widetilde{E}^{-1}(-\eta^+_{j})\right) D_{r}\right\}  \notag\\
&&+\int_{\zeta_{k(t)}}^{t}e^{A(t-s)}f(s)ds.\notag 
\end{eqnarray}
Also, if
\begin{eqnarray*}
&&\eta=\eta^+_k=\eta^-_k,\quad k\in\mathbb{Z},\quad \widehat{E}=(I+C) \widetilde{E}(\eta)\widetilde{E}^{-1}(-\eta),
\end{eqnarray*}
the solution of \eqref{sistema_idepcag_lineal_homogeneo_constante} is
$$w(t)=\widehat{W}(t,\tau)w(\tau),$$
where 
\begin{eqnarray*}
\widehat{W}(t,\tau)&=&\widetilde{E}(\eta(t))\widetilde{E}^{-1}(-\eta^+_{k(t)})\widehat{E}^{k(t)-k(\tau)-1}(I+C) \widetilde{E}(\eta)\widetilde{E}^{-1}(\eta(\tau)),
\end{eqnarray*}
is the fundamental matrix for \eqref{sistema_idepcag_lineal_homogeneo_constante} with $t\in I_{k(t)}$ and $\tau \in I_{k(\tau) }.$\\
The solution for \eqref{sistema_idepcag_lineal_no_homogeneo_constante} is given by
\begin{eqnarray}
y(t) &=&\widetilde{E}(\eta(t))\widetilde{E}^{-1}(-\eta^+_{k(t)})\widehat{E}^{k(t)-k(\tau)-1}( I+C) \widetilde{E}(\eta^-_{k(\tau)+1})\widetilde{E}^{-1}(\eta(\tau))\notag\\
&&\cdot \left(y(\tau)+\int_{\tau}^{\zeta_{k(\tau) }} e^{A(\tau-s)}f(s)ds\right)  \label{FORMULA_VARIACION:PARAMETROS_IDEPCAG_GENERAL_sistema_constante_2} \\
&&+\widetilde{E}(\eta(t))\widetilde{E}^{-1}(-\eta^+_{k(t)})\cdot \left\{\sum_{r=k(\tau)+1}^{k(t)}\widehat{E}^{k(t)-r} \int_{t_{r}}^{\zeta_{r}}e^{A(t_{r}-s)}f(s)ds\right.\notag\\
&&+\sum_{r=k(\tau)}^{k(t)-1}\widehat{E}^{k(t)-r} \int_{\zeta_{r}}^{t_{r+1}}(I+C)e^{A(t_{r+1}-s)}f(s)ds+\left.\sum_{r=k(\tau)+1}^{k(t)}\widehat{E}^{k(t)-r} D_{r}\right\}  \notag\\
&&+\int_{\zeta_{k(t)}}^{t}e^{A(t-s)}f(s)ds.\notag 
\end{eqnarray}
\end{enumerate}
\bigskip
\begin{remark}
\begin{enumerate}
\item[] 
\item We recover the variation of parameters concluded in \cite{P2011} when $D_{r}=C_{r}=0.$ 
\item Also, our result implies the variation of constant formulas given in section \ref{section_overview}
\end{enumerate}
\end{remark}

\section{Some Examples of Linear IDEPCAG systems}
In \cite{oztepe2012convergence}, H. Bereketoglu and G. Oztepe
studied the following linear IDEPCAG \begin{equation}
\begin{tabular}{ll}
$z^{\prime }(t)=A(t)\left(z(t)-z(\gamma (t))\right) +f(t),$ & $t\neq t_k$ \\ 
$z(t_{k})=z(t_{k}^{-})+D_k,$ & $t=t_k.$ \\ 
$z(\tau) =z_{0}$ & 
\end{tabular}
\label{SISTEMA_BEREK}
\end{equation}
where $\gamma(t)$ is some piecewise constant argument of generalized type,  $A(t)$ is a continuous locally integrable matrix,
$D:\mathbb{N}\rightarrow \mathbb{R}$ 
is such that $D_k\neq 0, \forall k\in \mathbb{N}.$
The authors originally considered the cases $\gamma_1(t)=[t+1],$ and $\gamma_2(t)=[t-1].$ Hence, $t_{k}=k, \zeta_{1,k}k=k+1$ and  $\zeta_{2,k}=k-1$, respectively.\\
Let $\Phi(t) $ be the fundamental matrix of the ordinary differential system
\begin{equation}
x^{\prime }(t)=A(t)x(t).  \label{sistema_ordinario_berek}
\end{equation}
It is well known that $\Phi ^{-1}(t) $ is the fundamental matrix of the adjoint system associated with \eqref{sistema_ordinario_berek}. So, it satisfies 
\begin{equation*}
\left(\Phi ^{-1}\right) ^{^{\prime }}(t)=-\Phi ^{-1}(t)A(t).
\end{equation*}
Therefore, we have
\begin{align*}
J(t,t_k)& =I-\int_{t_k}^{t}\Phi(t_k,s)A(s)ds \\
& =I+\Phi(t_k)\left(\int_{t_k}^{t}-\Phi^{-1}(s)A(s)ds\right) \\
& =I+\Phi(t_k)\left(\Phi^{-1}(t)-\Phi^{-1}(t_k)\right) \\
& =\Phi(t_k,t) \\
& =\Phi^{-1}(t,t_k),
\end{align*}
$E(t,t_k)=\Phi(t,t_k)J(t,t_k)=\Phi(t,t_k)\Phi^{-1}(t,t_k)=I,$ and, as a result of last estimations,  for $t,t^{\prime}\in I_{k},$ we have $W\left( t,t^{\prime }\right)=I$. Hence, the linear homogeneous IDEPCAG (is a DEPCAG because $C_k=0$) system
\begin{equation}
\begin{tabular}{ll}
$w^{\prime }(t)=A(t)\left(w(t)-w(\gamma (t))\right),$ & $t\neq t_k$ \\ 
$w(t_{k})=w(t_{k}^{-})$ & $t=t_k.$ \\ 
$w(\tau)=w_{0},$ & 
\end{tabular}
\label{SISTEMA_BEREK_homogeneo}
\end{equation}
has the constant solution $w(t)=w(\tau).$\\
Finally, for the variation of parameters formula \eqref{FORMULA_VARIACION_PARAMETROS_IDEPCAG_GENERAL_FINAL}, the solution for \eqref{SISTEMA_BEREK} is
\begin{eqnarray}
y(t) &=&y(\tau)+\int_{\tau}^{\zeta_{k(\tau) }}\Phi(\tau ,s)f(s)ds  +\sum_{r=k(\tau) +1}^{k(t)}\int_{t_{r}}^{\zeta_{r}}\Phi (t_{r},s)f(s)ds \notag \\
&&+\sum_{r=k(\tau)}^{k(t)-1}\int_{\zeta_{r}}^{t_{r+1}}\Phi(t_{r+1},s)f(s)ds+\int_{\zeta_{k(t)}}^{t}\Phi(t,s)f(s)ds+\sum_{r=k(\tau)+1}^{k(t)} D_{r},  \notag
\end{eqnarray}
\begin{figure}[h!]
\centering
\includegraphics[scale=0.3]{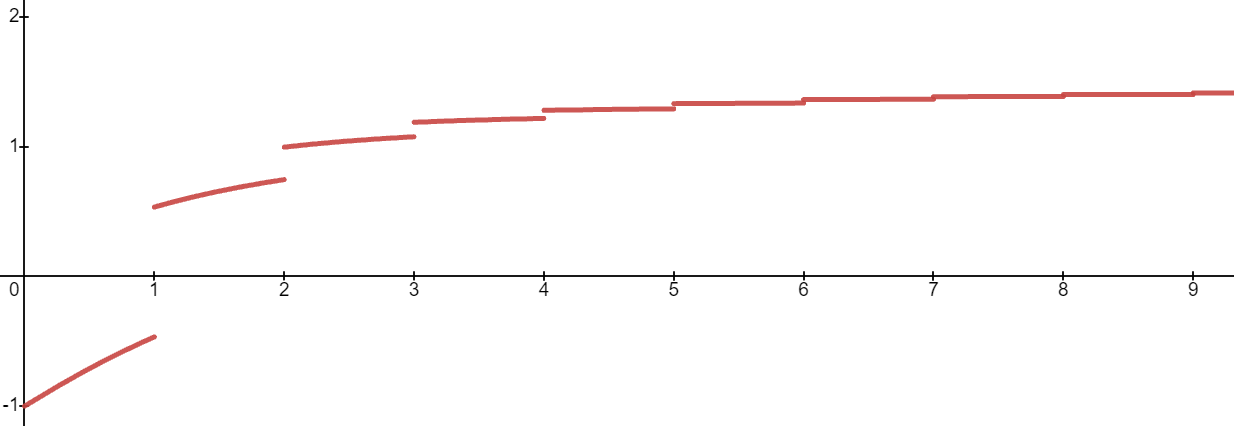}
\caption{Solution of  \eqref{SISTEMA_BEREK} with $\gamma(t)=[t]+7/10$, $D_r=1/r^2$, $A(t)=1/(t+1),$ $f(t)=\exp(-t)$ and $y(0)=y_0=-1.$}
\end{figure}
\begin{remark}
    This is the IDEPCAG case for the well-known differential equation studied by K.L Cooke and J.A. Yorke in \cite{COOKE_Yorke}. The authors investigated
    the following delay differential equation (DDE):
    \begin{equation*}
        x'(t)=g(x(t))-g(x(t-L)),
    \end{equation*}
    where $x(t)$ denotes the number of individuals in a population, the number of births is $g(x(t))$, and $L$ is the constant life span of the individuals in the population. Then, the number of deaths $g(x(t - L))$. Since the difference $g(x(t)) - g(x(t - L)]$ means the change of the population. Therefore $x'(t)$ corresponds to the growth of the population at instant $t$.\\
    In \eqref{SISTEMA_BEREK_homogeneo}, we considered $g(x(t))=A(t)x(t)$ and the constant delay in the Cooke-Yorke equation is regarded as a piecewise constant argument $\gamma(t)$. Notice that if $D_r$ is summable and $f(t)=0\,\,\forall t\geq\tau,$ then the solution of \eqref{SISTEMA_BEREK} tends to the constant $$y_\infty=y(\tau)+ \sum_{t_r\geq t_{k(\tau)+1}} D_{r}, \quad \text{ as } t\to\infty,$$
    no matter what $\gamma(t)$ was used.
    For further about asymptotics in IDEPCAG, see \cite{CTP2019}.
\begin{figure}[h!]
\centering
\includegraphics[scale=0.35]{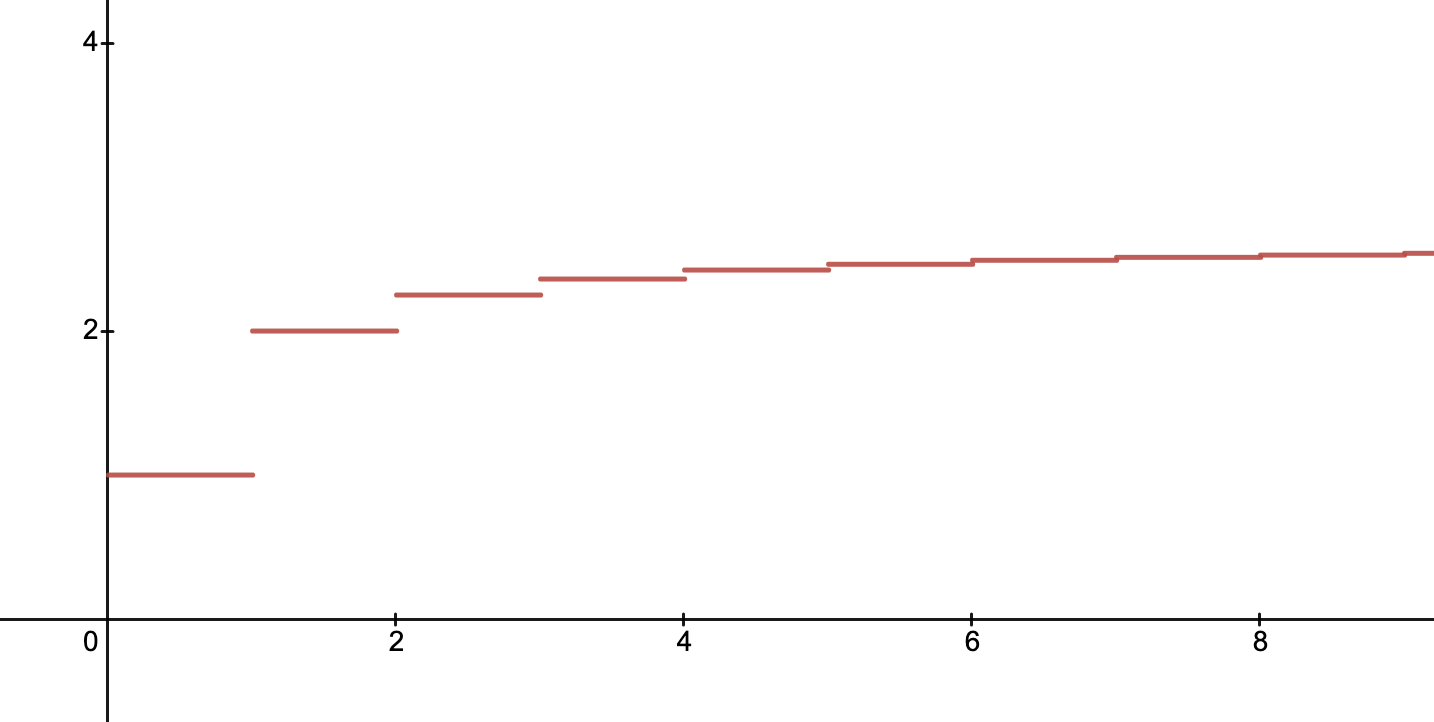}
\caption{Solution of \eqref{SISTEMA_BEREK} with $D_k=1/k^2,\,\,f(t)=0$ and $z(0)=1.$}
\end{figure}
\end{remark}

Let the following IDEPCA
\begin{equation}
\begin{tabular}{ll}
$z'(t)=\sin(2\pi t)z\left(\left[\frac{t}{h}\right]h+\beta h\right)+1,$ & $t\neq kh,\quad k\in\mathbb{N}$, \\ 
$z(kh)=\left(-\dfrac{1}{2}\right)z(kh^{-})+\dfrac{1}{2},$ & $t=kh,$ \\ 
$z(0)=z_{0}$. & 
\end{tabular}
\label{IDEPCA_EJEMPLO}
\end{equation}
where $h>0,\,0\leq \beta\leq 1.$\\
It is easy to see that $t_k=kh,\,\,\zeta_k=(k+\beta)h,\,\,k\in\mathbb{N}$ and 
$$I_k^+=[kh,(k+\beta)h],\quad I_k^-=[(k+\beta)h,(k+1)h).$$

We see that
\begin{eqnarray*}
\nu_k^{+}(\sin(2\pi t))&\leq&\beta h<1,\quad \text{if } h \text{ is small enough,}\\
\nu_k^{-}(\sin(2\pi t))&\leq& (1-\beta) h<1 ,\quad \text{if } h \text{ is small enough,}\\
E(t,\tau)&=&1+\int_{\tau}^{t}\sin(2\pi s)ds.\\
\label{MATRIZ_J_ejemplo}
\end{eqnarray*}
The fundamental matrix of the homogeneous equation associated with \eqref{IDEPCA_EJEMPLO} is
\begin{eqnarray*}
W(t,0)&=&\left(1+\int_{[t/h]h+\beta h}^{t}\sin\left(2\pi s\right)ds\right)\left(1+\int_{[t/h]h+\beta h}^{[t/h]h}\sin\left(2\pi s\right)ds\right)^{\left(-1\right)}\\
&&\cdot \left(-\frac{1}{2}\right)^{[t/h]}\left(\prod_{j=0}^{[t/h]-1}\left(1+\int_{(j+\beta)h}^{\left(j+1\right)h}\sin\left(2\pi s\right)ds\right)\left(1+\int_{(j+\beta)h}^{jh}\sin\left(2\pi s\right)ds\right)^{\left(-1\right)}\right).
\label{MATRIZ_FUNDAMENTAL_IDEPCA_ejemplo}
\end{eqnarray*}
Hence, the solution of \eqref{IDEPCA_EJEMPLO} is
\begin{eqnarray*}
z(t)&=&W(t,0)z_0+\left(-\frac{1}{2}\right)\left(1-\beta\right)h\sum_{r=0}^{[t/h]-1}W\left(t,\left(r+1\right)h\right)+\left(t-([t/h]h+\beta h)\right)\\
&&+W\left(t,0\right)\beta h+\beta h\sum_{r=1}^{[t/h]}W\left(t,rh\right)+
\left(-\frac{1}{2}\right)\sum_{r=0}^{[t/h]-1}W\left(t,\left(r+1\right)h\right).
\end{eqnarray*}
The piecewise constant used in this example was introduced in \cite{Torres_proyecciones} to study the approximation of solutions of differential equations (under some stability assumptions and taking $h\to 0$.)\\
\begin{figure}[h!]
\centering
\includegraphics[scale=0.25]{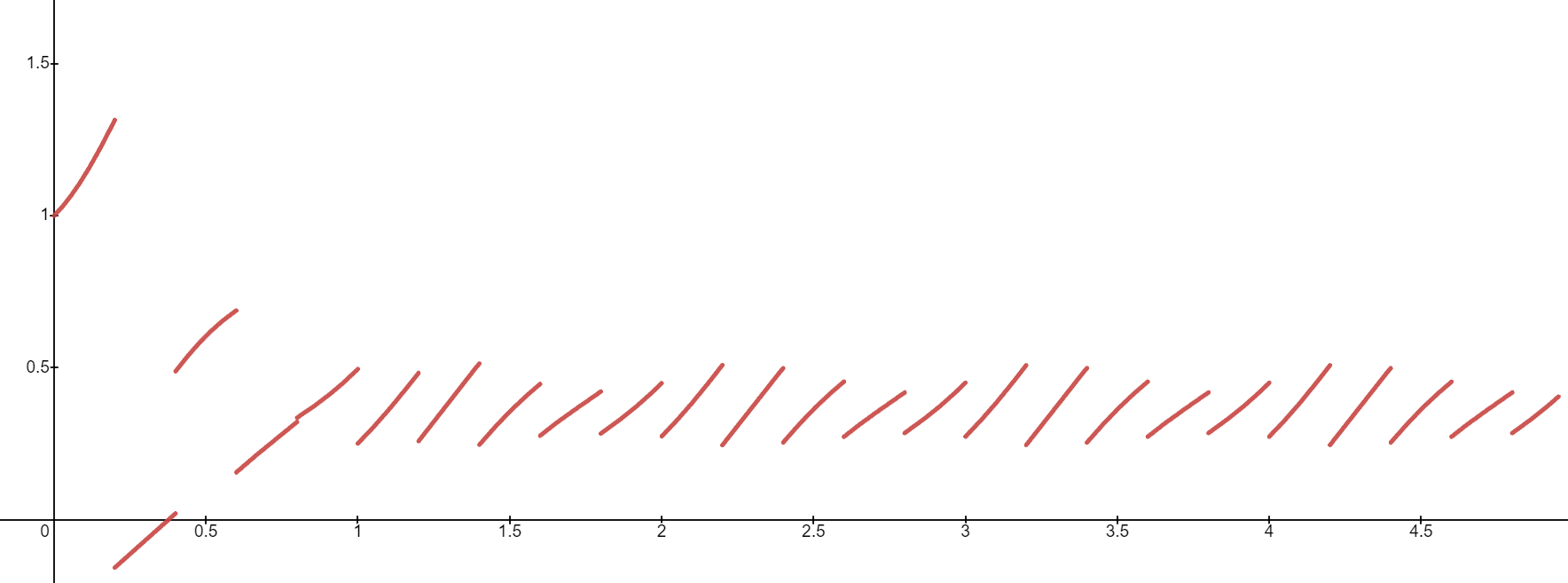}
\caption{Solution of  \eqref{IDEPCA_EJEMPLO} with $h=\beta=0,2$.}
\end{figure}
\section{Conclusions}
In this work, we gave a variation of parameters formula for impulsive differential equations with piecewise constant arguments. We analyzed the constant coefficients case and gave several examples of formulas applied to some concrete piecewise constant arguments. We extended some cases treated before and showed the effect of the impulses in the dynamic. 

\section*{Acknowledgments}
\noindent \emph{Manuel Pinto} thanks for the support of Fondecyt project 1170466. \\
\emph{Ricardo Torres} thanks to DESMOS PBC for granting permission to use the images employed in this work. They were created with the DESMOS graphic calculator \\ \url{https://www.desmos.com/calculator}.
\subsection*{Author contributions}
All of the authors contributed equally to this work.

\subsection*{Financial disclosure}
None reported.

\subsection*{Conflict of interest}
None declared.

\bibliographystyle{abbrv}
\bibliography{sca}

\begin{thebibliography}{10}

\bibitem{A04}
M.~Akhmet.
\newblock Stability of differential equations with piecewise constant arguments of generalized type.
\newblock {\em Nonlinear Analysis: Theory, Methods \& Applications}, 68(4):794--803, 2008.

\bibitem{AK3}
M.~Akhmet.
\newblock {\em Principles of Discontinuous Dynamical Systems}.
\newblock Springer, New York, Dordrecht, Heidelberg, London, 2010.

\bibitem{AK2}
M.~Akhmet.
\newblock {\em Nonlinear Hybrid Continuous-Discrete-Time Models}.
\newblock Atlantis Press, Amsterdam-Paris, 2011.

\bibitem{CTP2019}
S.~Castillo, M.~Pinto, and R.~Torres.
\newblock Asymptotic formulae for solutions to impulsive differential equations with piecewise constant argument of generalized type.
\newblock {\em Electron. J. Differential Equations}, 2019(40):1--22, 2019.

\bibitem{kuo2022}
K.-S. Chiu.
\newblock Periodic solutions of impulsive differential equations with piecewise alternately advanced and retarded argument of generalized type.
\newblock {\em Rocky Mountain Journal of Mathematics}, 52(1):87 -- 103, 2022.

\bibitem{Kuo2023}
K.-S. Chiu and I.~Berna.
\newblock Nonautonomous impulsive differential equations of alternately advanced and retarded type.
\newblock {\em Filomat}, 37:7813--7829, 2023.

\bibitem{COOKE1}
K.~Cooke and S.~Busenberg.
\newblock Models of vertically transmitted diseases with sequential-continuous dynamics.
\newblock In {\em Nonlinear Phenomena in Mathematical Sciences: Proceedings of an International Conference on Nonlinear Phenomena in Mathematical Sciences, Held at the University of Texas at Arlington, , June 16-20, 1980}, pages 179--187, New York, 1982. Academic Press.

\bibitem{CoWi87}
K.~Cooke and J.~Wiener.
\newblock An equation alternately of retarded and advanced type.
\newblock In {\em Proceedings of of the American mathematical society, April, 1987}, pages 726--732, New York, 1987.

\bibitem{CoWi84}
K.~L. Cooke and J.~Wiener.
\newblock Retarded differential equations with piecewise constant delays.
\newblock {\em J. of Mathematical Analysis and Applications}, 99(1):265--297, 1984.

\bibitem{COOKE_Yorke}
K.~L. Cooke and J.~A. Yorke.
\newblock Some equations modelling growth processes and gonorrhea epidemics.
\newblock {\em Mathematical Biosciences}, 16(1):75--101, 1973.

\bibitem{DAI}
L.~Dai.
\newblock {\em Nonlinear Dynamics of Piecewise Constant Systems and Implementation of Piecewise Constant Arguments}.
\newblock World Scientific Press Publishing Co, New York, 2008.

\bibitem{DAI_SINGH}
L.~Dai and M.~Singh.
\newblock On oscillatory motion of spring-mass systems subjected to piecewise constant forces.
\newblock {\em Journal of Sound and Vibration}, 173(2):217--231, 1994.

\bibitem{JaDeo}
K.~Jayasree and S.~Deo.
\newblock Variation of parameters formula for the equation of cooke and wiener.
\newblock In {\em Proceedings of the American Mathematical Society}, volume 112, United States, 1991.

\bibitem{MenYan}
Q.~Meng and J.~Yan.
\newblock Nonautonomous differential equations of alternately retarded and advanced type.
\newblock {\em International Journal of Mathematics and Mathematical Sciences}, 26:597--603, 2001.

\bibitem{203}
A.~Myshkis.
\newblock On certain problems in the theory of differential equations with deviating argument.
\newblock {\em Russian Mathematical Surveys}, 32(2):173--203, 1977.

\bibitem{oztepe2012convergence}
G.~Oztepe and H.~Bereketoglu.
\newblock Convergence in an impulsive advanced differential equations with piecewise constant argument.
\newblock {\em Bulletin of Mathematical Analysis and Applications}, 4(2012):57--70, 2012.

\bibitem{P2011}
M.~Pinto.
\newblock Cauchy and {G}reen matrices type and stability in alternately advanced and delayed differential systems.
\newblock {\em Journal of Difference Equations and Applications}, 17(2):235--254, 2011.

\bibitem{Samoilenko}
A.~Samoilenko.
\newblock {\em Impulsive Differential Equations}.
\newblock World Scientific Press, Singapur, 1995.

\bibitem{ShWi83}
S.~M. Shah and J.~Wiener.
\newblock Advanced differential equations with piecewise constant argument deviations.
\newblock {\em International J. of Mathematics and Mathematical Sciences}, 6(4):71--703, 1983.

\bibitem{Torres1}
R.~Torres.
\newblock Ecuaciones diferenciales con argumento constante a trozos del tipo generalizado con impulso.
\newblock Msc. thesis, Facultad de Ciencias, Universidad de Chile, Santiago, Chile, 2015.
\newblock Available at \url{https://repositorio.uchile.cl/handle/2250/188898}.

\bibitem{Torres_proyecciones}
R.~Torres, S.~Castillo, and M.~Pinto.
\newblock How to draw the graphs of the exponential, logistic, and gaussian functions with pencil and ruler in an accurate way.
\newblock {\em Proyecciones, Journal of Mathematics. (to appear)}, 2023.

\bibitem{WiAf88}
J.~Wiener and A.~R. Aftabizadeh.
\newblock Differential equations alternately of retarded and advanced type.
\newblock {\em Journal of Mathematical Analysis and Applications}, 1(129):243--255, 1988.

\end{thebibliography}

\end{document}